\documentclass[reqno,12pt]{amsart}
\usepackage{amssymb, amsmath,latexsym,amsfonts,amsbsy, amsthm}
\usepackage{xcolor}
\usepackage{txfonts}
\usepackage{hyperref}
\usepackage{graphicx}

\allowdisplaybreaks[4]

\makeatletter
\renewcommand\normalsize{%
  \@setfontsize\normalsize{13.5pt}{17pt}%
  \abovedisplayskip 10pt plus 2pt minus 5pt
  \abovedisplayshortskip \z@ plus 3pt
  \belowdisplayshortskip 6pt plus 3pt minus 3pt
  \belowdisplayskip \abovedisplayskip
  \let\@listi\@listI}
\makeatother

\newtheorem{theorem}{Theorem}[section]
\newtheorem{definition}[theorem]{Definition}
\newtheorem{lemma}[theorem]{Lemma}
\newtheorem{proposition}[theorem]{Proposition}

\newtheorem{remark}[theorem]{Remark}
\newtheorem{problem}[theorem]{Problem}

\def\bR{\mathbb R}

\def\bZ{\mathbb Z}

\def\ve{\varepsilon}
\def\la{\lambda}

\def\t{\tilde}
\def\q{\quad}
\def\qq{\qquad}
\def\th{\theta}

\def\dl{\delta}
\def\Dl{\Delta}
\def\lt{\left}
\def\les{\lesssim}
\def\rt{\right}

\def\i{\infty}
\def\e{\varepsilon}

\def \ls{\lesssim}
\def\p{\partial}
\def\f{\frac}
\def\na{\nabla}
\def\al{\alpha}

\def\O{\Omega}
\def\o{\omega}

\def\s{\sqrt}
\def\nn{\nonumber}
\def\be{\begin{equation}}
\def\ee{\end{equation}}
\def\bes{\begin{equation*}}
\def\ees{\end{equation*}}
\def\bali{\begin{aligned}}
\def\eali{\end{aligned}}
\def\bl{\boldsymbol}

\def\pr{\prime}
\def\ed{\buildrel\hbox{\footnotesize def}\over =}

\def\ba{\begin{equation}\begin{aligned}}
	\def\ea{\end{aligned}\end{equation}}
\def\bn{\[\begin{aligned}}
\def\en{\end{aligned}\]}
\def\l{\label}
\def\ed{\buildrel\hbox{\footnotesize def}\over =}

\def\t{\tilde}
\def\cd{\cdot}
\def\mL{\mathcal{L}}
\def\mN{\mathcal{N}}
\def\mT{\mathcal{T}}
\def\mf{\mathrm{f}}
\def\s{\sigma}
\def\un{\underbrace}
\def\bm{\boldsymbol}

\newcommand{\man}[1]{%
\begin{align}#1\end{align}%
}

\begin{document}

\title[Wall law in a wedge]
{A wall-law approximation for Jeffery-Hamel flows in rough wedges}
\author[Z. Li]{Zijin Li}
\address[Z. Li]{School of Mathematics and Statistics, Nanjing University of Information Science and Technology, Nanjing 210044, China}
\email{zijinli@nuist.edu.cn}

\author[X. Pan]{Xinghong Pan}
\address[X. Pan]{School of Mathematics and Key Laboratory of MIIT, Nanjing University of Aeronautics and Astronautics, Nanjing 211106, China}
\email{xinghong\_87@nuaa.edu.cn}

\begin{abstract}
 We study the two-dimensional stationary Navier-Stokes equations in an unbounded wedge with small rough perturbations of its angular boundaries. The Jeffery-Hamel flow in the corresponding straight wedge is taken as the effective background flow. Under a suitable small-flux condition,  we prove the existence of weak solutions and establish an $H^1$-type energy error estimate of order $O(\e^2)$. For sufficiently small wedge angles, we further derive weighted estimates and improve the squared $L^2$-error to $O(\e^3)$.
\end{abstract}



\maketitle

\setcounter{equation}{0}
\numberwithin{equation}{section}

\section{Introduction}

 Flows in domains with rough boundaries arise naturally in many problems of fluid mechanics, where the microscopic geometry of the wall may have a non-negligible influence on the macroscopic behavior of the flow. Understanding how boundary roughness modifies the effective motion of viscous incompressible fluids is a classical issue closely related to wall laws, homogenization, and boundary layer analysis. A wall law replaces a genuinely rough boundary by an effective boundary condition imposed on a smooth artificial interface, thereby encoding the averaged influence of microscopic wall geometry in a macroscopic and computationally tractable form. Typical effective conditions include corrected no-slip laws and Navier-type slip laws. Since resolving the full rough geometry is often prohibitively expensive, wall laws provide reduced models that preserve the leading effects of roughness on the velocity field, boundary-layer structure, and drag. From the analytical point of view, an important problem is to justify these effective laws rigorously and to quantify the associated error estimates.

Effective boundary conditions for viscous flows over rough walls have been studied from asymptotic, numerical, and rigorous viewpoints. Early work of Achdou, Pironneau and Valentin \cite{APV1998JCP} proposed first-order and second-order wall laws for laminar flows over periodic rough boundaries. Rigorous derivations of Navier-type effective conditions were later obtained by J\"ager and Mikeli\`c \cite{JM2001JDE,JM2003CMP} for incompressible viscous flows over periodic rough surfaces and for Couette configurations with riblets. Random roughness was treated by Basson and G\'{e}rard-Varet \cite{BG2008CPAM}, who derived wall laws for stationary Navier--Stokes flows in randomly rough channels, and by G\'{e}rard-Varet \cite{G2009CMP}, who established quantitative error estimates. Higher-order multiscale wall laws were investigated by Bresch and Milisic \cite{BM2010Quart}, while Bucur, Feireisl, and Ne\v casov\'a \cite{BFN2010ARMA} clarified the relationship between wall roughness and effective slip behavior through a $\Gamma$-convergence approach. The relevance and limitations of slip laws near arbitrary irregular boundaries, as well as effective boundary conditions starting from microscopic slip, were further studied in \cite{GVM2010CMP,DGV2011JDE}. Nonstationary extensions can be found, for instance, in \cite{Higaki2016JDE}.

Most of the above works concern flat boundaries, channels, or locally periodic roughness patterns. The present paper considers a different geometric setting: an unbounded wedge domain with rough angular boundaries and a flux condition at the apex. The limiting background flow is the Jeffery--Hamel flow, a classical radial solution in a wedge, whose singular spatial scaling may interact nontrivially with the rough boundary perturbation. Moreover, the roughness considered here is not a uniformly $C^1$ perturbation of a flat graph. After the logarithmic polar transformation, the angular perturbation becomes exponentially localized in the longitudinal variable. This feature is crucial in our estimates and leads to improved powers of the roughness parameter. Our goal in this paper is to establish wall-law type estimates for perturbations of Jeffery--Hamel flows in rough wedge domains. The analysis combines the logarithmic formulation, divergence corrections of Bogovskii type, compactness and fixed-point arguments, weighted energy estimates, and an adjoint method for sharper $L^2$ bounds.

Consider the stationary Navier--Stokes equations in a two-dimensional wedge domain:
\be\label{nschannel}
\lt\{
\bali
&\bl{u}^\e\cdot\na \bl{u}^\e+ \na \Pi^\e- \Dl \bl{u}^\e=0, \q x\in \t{\O}^\e,\\
&\text{div}\,\bl{u}^\e=0,\q \q x\in \t{\O}^\e,
\eali
\rt.
\ee
subject to the boundary condition
\ba\l{BC}
\bl{u}=0\q \text{on}\q \p\t{\O}^\e,
\ea
where
\[
\t{\O}^\e\ed \lt\{x=(x_1,x_2)\in \bR_+\times\bR\,:\,-\al - \t{\th}^\e_-(r) <\th<\al + \t{\th}^\e_+(r) \rt\}\,,
\]
\begin{figure}[!ht]
	\centering
	\includegraphics[scale=.35]{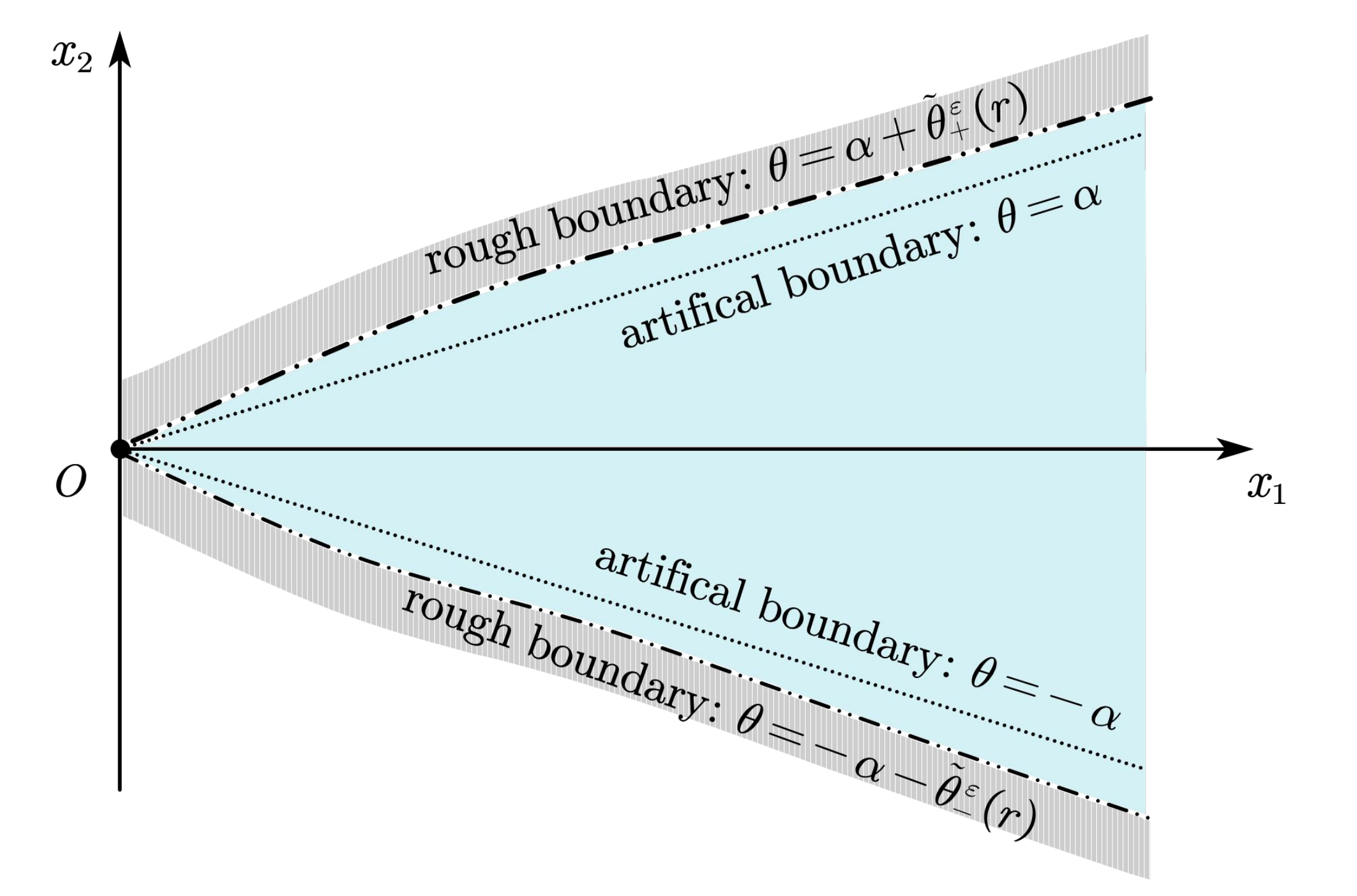}
	\caption{The wedge domain $\t{\O}^\e$ in $\mathbb{R}^2$}
	\label{Wedge}
\end{figure}
with half angle $\al\in (0,\f{\pi}{2})$. Here and below
\[
r=\sqrt{x_1^2+x_2^2}\,,\q\th=\arctan\f{x_2}{x_1}
\]
denote the radial distance and the polar angle, respectively. The upper and lower boundaries are perturbed by rough oscillations:
 \be\label{scale}
 \t{\th}^\e_\pm(r)\ed\lt\{
 \bali
 &\e\t{\gamma}_{\pm}(r^\f{1}{\e}), \hskip .2cm\q \text{for } r\geq 1,\\
 &\e\t{\gamma}_{\pm}(r^{-\f{1}{\e}}), \q \text{for } r< 1,
 \eali
 \rt.
 \ee
where
\ba\l{GamSet}
\t{\gamma}_\pm(\rho)\in C^1[1,\infty)
\ea
denote the upper and lower rough perturbation profiles, respectively. Moreover, they decay algebraically to zero at infinity, i.e.
\be\l{BSet}
\t{\gamma}_\pm(\rho) = O(\rho^{-\kappa}),\q \text{ as }\q \rho\rightarrow \i,
\ee
for some $\kappa>0$. Intuitively, the rough boundaries $\t{\th}_{\pm}^\e$ converge to the corresponding smooth effective boundaries, referred to as wall-law boundaries $\{\th=\pm\al\}$ as $\e\to 0$.

Compared with the standard oscillatory graph perturbations of the form $x_2=\e \gamma\left(\f{x}{\e}\right)$ considered, for instance, in \cite{BG2008CPAM,G2009CMP,BM2010Quart}, the angular perturbation in \eqref{scale} needs not remain uniformly bounded in $C^1$ unless an appropriate decay assumption is imposed on $\gamma^\prime$. Indeed, if $\gamma^\prime$ is only bounded,
\bes
(\t{\th}^\e_\pm)^\prime(r)=\lt\{
\bali
&\t{\gamma}_\pm^\prime(r^\f{1}{\e})r^{\f{1}{\e}-1}\rightarrow \i, \q\hskip .9cm\text{ for fixed $r>1$, as $\e\rightarrow 0$},\\
&-\t{\gamma}_\pm^\prime(r^{-\f{1}{\e}})r^{-\f{1}{\e}-1}\rightarrow \i, \q\text{ for fixed $r<1$, as $\e\rightarrow 0$}.
\eali
\rt.
\ees
Thus the present roughness model is adapted to wedge geometry and allows boundary oscillations that are rougher than uniformly $C^1$ graph perturbations.

Finally, we supplement\eqref{nschannel}--\eqref{BC} with the prescribed-flux condition:
\ba\l{flux}
\int_{\Sigma}\bl{u}\cdot\bl{n}\,dS=\Phi.
\ea
Here $\Sigma$ is a $C^1$ simple curve in $\t{\O}^\e$ connecting the upper and lower boundaries, $\bl{n}$ is the unit normal vector to this curve, pointing away from the origin. By the divergence theorem and the incompressibility of the fluid, $\Phi$ is independent of the choice of $\Sigma$. The flux can be interpreted as the strength of the fluid source (or sink if $\Phi<0$) at the apex of the wedge. In this paper we restrict attention to the source-flow case $\Phi>0$ since the opposite case is similar.

\subsection{Reformulation of the problem in the polar coordinates}

Since the reference domain is a perturbed sector, it is more convenient to reformulate our system in polar coordinates. Denote the velocity by
\bn
\bl{u}^\e=\t{u}^\ve(r,\th) \bl{e}_r+\t{v}^\e(r,\th) \bl{e}_\th\,,
\en
with
\[
\bl{e}_r\ed(\cos\th,\sin\th)^{\mathrm T},\q\bl{e}_\th\ed(-\sin\th,\cos\th)^{\mathrm T},
\]
we reformulate \eqref{nschannel}--\eqref{BC}--\eqref{flux} as follows:

\be\label{nspolar}
\left\{
\begin{aligned}
&(r \t{u}^{\e}\p_r+\t{v}^{\e}\p_\th) \t{u}^{\e}-\left(v^{\e}\right)^2+r \p_r {\Pi}^{\e}-\left(\frac{\p^2_\th \t{u}^{\e}}{r}+\p_r(r\p_r \t{u}^{\e})-\frac{2}{r} \p_\th \t{v}^{\e}-\frac{\t{u}^{\e} }{r}\right)=0, \text{ in } \t{\O}^\e,\\
&(r \t{u}^{\e}\p_r+\t{v}^{\e}\p_\th) \t{v}^{\e} +\t{u}^{\e} \t{v}^{\e}+\p_\th {\Pi}^{\e}-\left(\frac{\p^2_\th \t{v}^{\e}}{r}+\p_r(r\p_r \t{v}^{\e})+\frac{2}{r} \p_\th \t{u}^{\e}-\frac{\t{v}^{\e}}{r}\right)=0, \hskip .3cm\text{ in } \t{\O}^\e,\\
&\p_\th \t{v}^{\e}+\p_r\left(r \t{u}^{\e}\right)=0,  \q \text{in}\q \t{\O}^\e,\\
&\lt(\t{u}^{\e},\t{v}^\e\rt)=0, \hskip 1.1cm\q \text{on}\q \p\t{\O}^\e.
\end{aligned}
\right.
\ee

For convenience of the subsequent calculations, we introduce the following change of variables: Let
\ba\l{Trans}
s\ed\log r,
\ea
and denote
\bn
u(s,\th)\ed\t{u}(r,\th)\,,\q v(s,\th)\ed\t{v}(r,\th)\,,\q p(s,\th)\ed\Pi(r,\th)\,.
\en
In this $(s,\th)$-coordinates, the following identities hold:
\be\l{STHCO}
r\p_r= \p_s,\q \Dl=\p^2_r+\f{1}{r}\p_r+\f{1}{r^2}\p^2_\th=e^{-2s}(\p^2_s+\p^2_\th)\ed e^{-2s}\Dl_s.
\ee

In view of \eqref{STHCO}, we rewrite \eqref{nspolar} as
\be\label{nslog}
\left\{
\begin{aligned}
&(u^{\e}\p_s+v^{\e}\p_\th) u^{\e} -\left(v^{\e}\right)^2+\p_s p^{\e}- e^{-s}\left(\Dl_s u^\e-2\p_\th v^{\e}-{u^{\e} }\right)=0, \q \text{in}\q \O^\e,\\
&(u^{\e}\p_s+v^{\e}\p_\th) v^{\e} +u^{\e} v^{\e}+\p_\th p^{\e}- e^{-s}\left(\Dl_s v^\e+2 \p_\th u^{\e}-{v^{\e}}\right)=0,\q \hskip .15cm\text{in}\q \O^\e,\\
&\p_s u^{\e}+u^\e+\p_\th v^\e=0,\q \text{in}\q \O^\e,\\
&\lt(u^{\e},v^\e\rt)=0, \q \text{on}\q \p\O^\e,\\
&\int^{\al+\th^\e_+(s)}_{-\al-\th^\e_{-}(s)} e^su^\e(\th) d\th=\Phi,\q \text{ for  any } s\in\bR,
\end{aligned}
\right.
\ee
where the domain is rewritten as
\[
\O^\e\ed \lt\{(s,\th)\in \bR^2\,:\,-\al -{\th}^\e_-(s) <\th<\al + {\th}^\e_+(s) \rt\}.
\]
\begin{figure}[!ht]
	\centering
	\includegraphics[scale=.3]{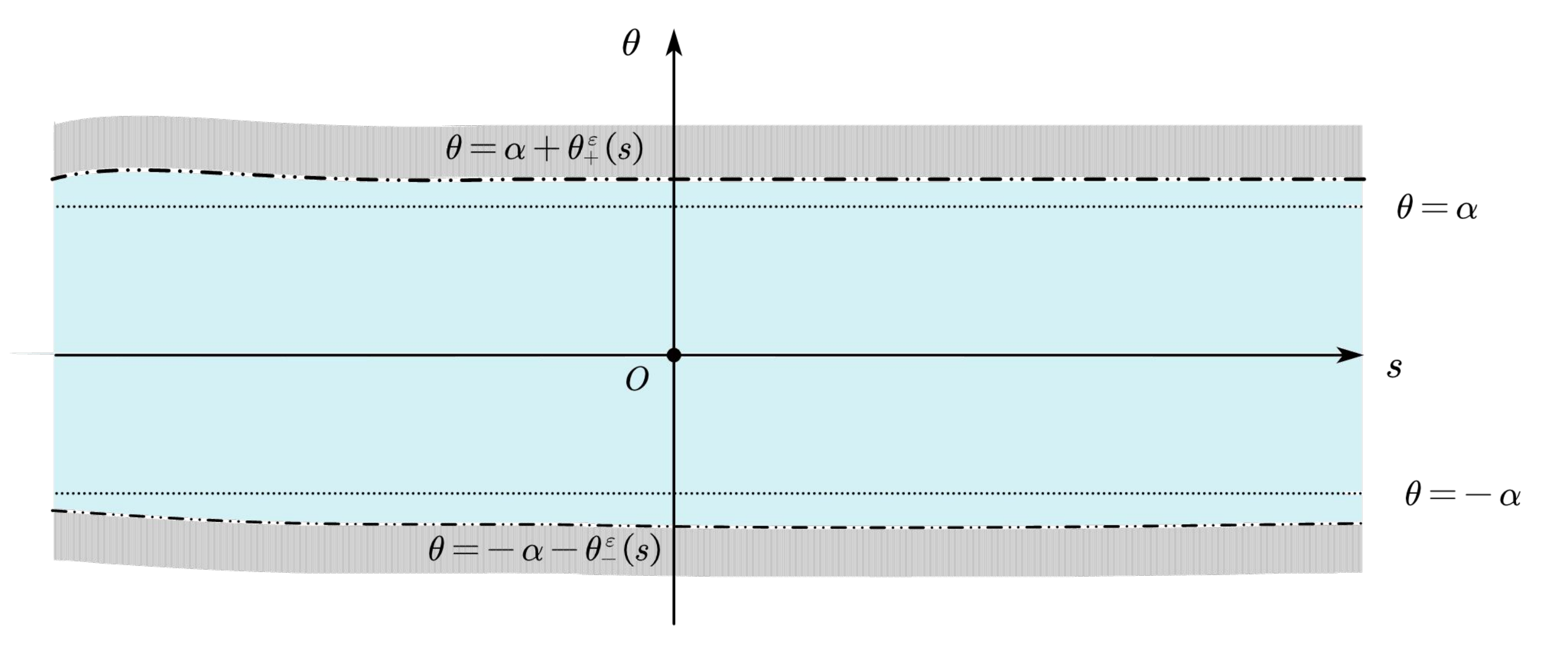}
	\caption{The wedge domain $\O^\e$ in $\mathbb{R}^2$}
	\label{Wedge1}
\end{figure}
Here, under the coordinate transformation \eqref{Trans}, the boundary profile becomes
\[
\gamma_\pm(s)\ed\lt\{
\bali
&\t{\gamma}_{\pm}(e^s), \hskip.2cm\q s\geq 0,\\
&\t{\gamma}_{\pm}(e^{-s}), \q s< 0.
\eali
\rt.
\]
We note that, in view of \eqref{GamSet}--\eqref{BSet}, $0\leq {\gamma}_\pm(s)\leq1$ and
\be\l{Ass}
{\gamma}_\pm(s) = O(e^{-\kappa |s|}),\q \text{ as } s\rightarrow \pm\i.
\ee
The corresponding rough boundaries are therefore given by
\ba\l{VETH}
{\th}^\e_\pm(s)=\t{\th}^\e_\pm(r)=\e{\gamma}_{\pm}\big(\f{s}{\e}\big)\,.
\ea

From \eqref{Ass}--\eqref{VETH}, a direct calculation gives the following asymptotic behavior of $\th_{\pm}^\e$:
\be\label{exponetialdecay}
\th_\pm^\e(s)\les\e e^{-\f{\kappa |s|}{\e}}\,,\q\text{as}\q s\to\pm\infty\,.
\ee

When $\e=0$, system \eqref{nslog} is written as
\be\label{nslog1}
\left\{
\begin{aligned}
&(u^{0}\p_s+v^{0}\p_\th) u^{0} -\left(v^{0}\right)^2+\p_s p^0- e^{-s}\left(\Dl_s u^0-2\p_\th v^{0}-{u^{0} }\right)=0, \q \text{in}\q \O,\\
&(u^{0}\p_s+v^{0}\p_\th) v^{0} +u^{0} v^{0}+\p_\th p^0- e^{-s}\left(\Dl_s v^0+2 \p_\th u^{0}-{v^{0}}\right)=0,\q \hskip.2cm\text{in}\q\O,\\
&\p_s u^{0}+u^0+\p_\th v^0=0,  \q \text{in}\q \O,\\
&(u^{0},v^0)=0, \q \text{on}\q \p\O,\\
&\int^{\al}_{-\al} u^0(\th)e^s d\th=\Phi,\q \text{ for  any } s\in\bR.
\end{aligned}
\right.
\ee
Here $\O\ed\O^0$. System \eqref{nslog1} has a well-known special solution: the \emph{Jeffery--Hamel flow}
\[
(u^0,v^0)=(e^{-s}f(\th), 0),
\]
named after Jeffery \cite{Jeffery} and Hamel \cite{Hamel}, who studied such flows in the early twentieth century. Here $f$ is a function of $\th$ that solves
\be\label{jefferyhamelsol}
\left\{
\begin{aligned}
&-e^{-2s}\lt(f^2(\th)+f^{\pr\pr}(\th)\rt)+\p_s p^0=0,\\
&\p_\th p^0-2e^{-2s}f'(\th)=0,\\
&f(\pm\al)=0,\q \int^\al_{-\al} f(\th) d\th=\Phi.
\end{aligned}
\right.
\ee
Decoupling $f$ and the pressure $p^0$ in \eqref{jefferyhamelsol}, we find that
\be\label{jefferyhamelsol1}
\left\{
\begin{aligned}
&f^{\pr\pr}+4 f+f^2=C_{\Phi,\al},\\
&f(\pm\al)=0,\q \int^\al_{-\al}f(\th) d\th=\Phi. \\
&p^0(s,\th)= 2 e^{-2s}f(\th)-\f{1}{2}e^{-2s} C_{\Phi,\al}.
\end{aligned}
\right.
\ee
The constant $C_{\Phi,\al}$ is determined by the fixed angular $\al$ and the flux $\Phi$.

\subsection{Main results}

Before presenting our main results, we introduce some notation to simplify the presentation. Throughout this paper, $C_{a,b,c,...}$ denotes a positive constant depending on $a,\,b,\, c,\,...$, which may be different from line to line. $A\ls_{a,b,c,\dots} B$ means $A\leq C_{a,b,c,...} B$. For a norm $\|\cdot\|$, we use $\|(f,g,\cdots)\|$ to express
$\|f\|+\|g\|+\cdots$. The gradient $\na$ in our setting denotes $(\p_s, \p_\th)$.
For $1\leq p\leq\infty$ and $k\in\mathbb{N}$, $L^p$ denotes the usual Lebesgue space with norm
\[
\|f\|_{L^p(D)}:=
\lt\{
\begin{aligned}
&\left(\int_{D}|f(x)|^pdx\right)^{1/p},\quad &1\leq p<\infty,\\[3mm]
&\mathrm{ess sup}_{x\in D}|f(x)|,\quad &p=\infty,\\
\end{aligned}
\rt.
\]
while $W^{k,p}$ denotes the usual Sobolev space with its norm
\[
\begin{split}
\|f\|_{W^{k,p}(D)}:=&\sum_{0\leq|L|\leq k}\|\nabla^L f\|_{L^p(D)},\\
\end{split}
\]
where $L=(l_1,l_2)$ is a multi-index. We also simply denote $W^{k,p}$ by $H^k$ provided $p=2$. For simplicity, we use $L^2_s$ and $L^2_\th$ to denote the $L^2$ norms in the $s$ and $\th$ variables.

We now state the main results of the paper. For convenience, in the statements of the next two theorems, $(u^0,v^0)$ also denotes its zero extension in $\O^\e/\O$.
\begin{theorem}\label{thmexistence1}
There exists a constant $C_\ast$, independent of $\al$, $\Phi$, $\kappa$, and $\e$, such that if $\e\leq \min\{\f{\pi-2\al}{4}, \al, \f{\kappa}{4}\}$, and $\al$, $\Phi$ satisfy
\bn
\al\Phi + \f{\Phi}{(\pi-2\al)^2}\leq C_\ast\,,
\en
the system \eqref{nslog} has a weak solution satisfying
\begin{align}
&\big\|\p_s \Big((u^\ve-u^0, v^\e\Big)\big\|_{L^2(\O^\e)}^2+\f{\lt(\pi-2\al\rt)^2}{32\pi^2}\big\|\p_\th \Big((u^\ve-u^0,v^\ve\Big)\big\|_{L^2(\O^\e)}^2\leq C \f{\Phi^2\e^2}{\al^4\kappa(\pi-2\al)^2}\,. \label{estimate1}
\end{align}
Here $C$ is an absolute constant, which is independent of $\al$, $\Phi$, $\kappa$, and $\e$.
\end{theorem}
\qed
\begin{remark}
Comparing estimate \eqref{estimate1} with the first estimate in \cite[Theorem 2.2, equation (2.3)]{BG2008CPAM}, we gain one additional order in
$\e$. The reason is that the error distances decay exponentially, as established in \eqref{exponetialdecay}, together with their square roots have $L^2$ bounds of order $\e$. See \eqref{FR1} and \eqref{FR2}. Consequently, our estimate improves upon (3.8) in \cite{BG2008CPAM} by a factor of order $\e^{\f{1}{2}}$.
\end{remark}

Notice that in Theorem \ref{thmexistence1}, for fixed $\al\in(0,\f{\pi}{2})$, $\Phi>0,\kappa>0$, it follows that the kinetic energy for the error velocity in the middle strip $\O$ is of order $\e^2$. Furthermore, for suitably small angle $\al$, we can construct suitable weak solutions of \eqref{nslog} with a higher-order kinetic energy estimate in the strip $\O$. The proof is based on a set of weighted energy estimates. In this case, the following decay assumption for the derivative of the perturbed boundary function
\be\label{decayassumption}
M_1:=\sup_{s\geq 1} s|\t{\gamma}'(s)|<+\infty
\ee
appears to be necessary to ensure the rough boundary is uniformly Lipschitz with $0<\e<<1$ in $(s,\th)-$coordinates. In this way the Bogovskii constant $\frak{C}_\ast$ in Lemma \ref{div} is independent of $\e$ and $\al$. The second main result is as follows:

\begin{theorem}\label{thmexistence2}
 Assume that \eqref{decayassumption} holds. There exists a suitably small positive constant ${C}_\ast$, independent of $\al$, $\Phi$, $\kappa$, and $\e$, such that if $\e\leq\min\{
    \al^3, \f{\kappa}{12}\}$, and
\bn
\al +\Phi\leq {C}_\ast,
\en
the system \eqref{nslog} has a weak solution satisfying the following weighted estimate
\begin{align}
&\big\|\na\Big( u^\e-u^0, v^\e\Big)\cosh(s)\big\|_{L^2(\O^\e)}^2+\f{1}{\al^2}\big\|\Big(u^\e-u^0, v^\e\Big)\cosh s\big\|_{L^2(\O^\e)}^2\leq C\al^{-4}\kappa^{-1}\Phi^2\e^2\,. \label{weightedesti}
\end{align}

Moreover, based on the weighted estimate in \eqref{weightedesti}, we derive the following refined estimate, which gains one additional order in $\ve$ compared with the corresponding estimate in \eqref{estimate1}:
\begin{align}
&\|(u^\e-u^0,v^\e)(s,\pm\al)\cosh s\|^2_{L^2_s}\leq C\al^{-4} \Phi^2 \kappa^{-1} \e^3, \label{weightedweak1y}\\
&\|(u^\e-u^0,v^\e)\|^2_{L^2(\O)}\leq C \al^{-6}\Phi^2\left(1+\kappa^{-2}\right)\e^3. \label{weightedweak1x}
\end{align}
Here $C$ is an absolute constant, which is independent of $\al$, $\Phi$, $\kappa$, and $\e$.
\end{theorem}
\qed

\begin{remark}
Under the assumption that the angle $\al$ is sufficiently small, the estimates \eqref{weightedweak1y} and \eqref{weightedweak1x} improve by one additional order in $\e$ compared with the second and third estimates in equation (2.3) of Theorem 2.2 in \cite{BG2008CPAM}. The smallness assumption on the angle is needed to absorb the non-vanishing pressure term within our weighted functional framework, which is essential for deriving \eqref{weightedweak1x} via the adjoint method. See Section \ref{SEC5} below.
\end{remark}

Before closing this introduction, we briefly comment on the case of random boundary
perturbations. Although this setting is not pursued in the present paper, the deterministic
framework developed here still applies with only minor modifications; see Remark \ref{RMK1} below.

\begin{remark}\l{RMK1}
The parameter in the present random setting is the realization $\omega\in\Omega$, not a
small deterministic quantity $\varepsilon$. Hence, for each fixed $\omega$, one obtains a
deterministic problem of the same form as above, with random boundary data $\gamma(\omega)$ or,
equivalently, on a random perturbed domain $D(\omega)$. Writing
\[
u(\omega)=\mathcal S(\gamma(\omega)),
\]
where $\mathcal S:X\to Y$ denotes the deterministic solution operator, the only additional
issue is the measurability of the map $\omega\mapsto u(\omega)$. A sufficient condition is
that $\gamma:\Omega\to X$ be measurable and that $\mathcal S:X\to Y$ be Borel measurable; this
holds in particular when $\mathcal S$ is continuous. Alternatively, one may argue at the
level of suitable approximate problems, prove measurability for the corresponding
approximating solutions, and then pass to the limit. This is essentially the type of
argument used in \cite{BG2008CPAM,G2009CMP}. We omit these standard measurability
considerations here.

Unlike random roughness models, the approach in \cite{BFN2010ARMA} is deterministic: the
boundary perturbation is represented by a family of domains $\Omega^\e$, and the
effective boundary condition is recovered by $\Gamma$-convergence methods. Another
deterministic line of research is pursued by Bresch and Milisic
\cite{BM2010Quart}, in the periodic setting. There, the rough boundary is not
modeled probabilistically, but through a periodic microstructure, and the corresponding
wall laws are obtained from boundary layer analysis and multiscale expansions. In
particular, their work emphasizes the role of oscillatory corrections and effective
boundary conditions on a flat interface.

\end{remark}

The rest of this paper is organized as follows. In Section \ref{SEC2}, we construct the Jeffery--Hamel flow for sufficiently small flux by means of the Banach fixed-point theorem. Section \ref{SEC3} is devoted to the weak formulation of the wall-law problem, whose solvability is established in Section \ref{SEC4}. In Section \ref{SEC5}, we further derive the existence of weak solutions in a weighted energy space and give sharper kinetic bounds for the perturbation when the wedge angle is sufficiently small.

\section{Properties of the Jeffery--Hamel flow under a small-flux condition}\l{SEC2}

In this section, we solve \eqref{jefferyhamelsol1}$_{1,2}$ for sufficiently small flux by using the Banach fixed-point theorem. Differentiating both sides of \eqref{jefferyhamelsol1}$_1$, we obtain
\be\label{jh1}
\left\{
\begin{aligned}
&f^{\pr\pr\pr}+4 f'+2f f'=0,\\
&f(\pm\al)=0,\q \int^\al_{-\al}f(\th) d\th=\Phi. \\
\end{aligned}
\right.
\ee
The main result for system \eqref{jh1} is as follows.
\begin{proposition}\label{jh1prop}
For $\al\in (0,\f{\pi}{2})$, there exists a constant $\frak{c}_0$, independent of $\alpha$ and $\Phi$, such that when $\al\Phi\leq \frak{c}_0$, the system \eqref{jh1} has a $C^1[-\alpha,+\alpha]$ solution satisfying
\be\label{jh2ex}
\al\|f\|_{L^\i}+ \al^2\|f^\pr\|_{L^\i}\leq C\Phi.
\ee
where the constant $C$ is independent of $\al$ and $\Phi$.
 \end{proposition}
For convenience, we make the following change of variables: set $g(\th)=\f{\al}{\Phi} f(\al\th)$, then $g(\th)$ satisfies
\be\label{jh2g}
\left\{
\begin{aligned}
&g^{\pr\pr\pr}+4\al^2 g^\pr+2\al\Phi g g^\pr=0,\\
&g(\pm1)=0,\q \int^1_{-1}g(\th) d\th=1. \\
\end{aligned}
\right.
\ee
Thus Proposition \ref{jh1prop} follows directly from the following lemma:
\begin{lemma}
For $\al\in (0,\f{\pi}{2})$, there exists a constant $\frak{c}_0>0$, independent of $\al$ and $\Phi$, such that when $\al\Phi\leq \frak{c}_0$, the system \eqref{jh2g} has a  $C^3[-1,1]$ solution satisfying
\be\label{jh2exg}
\|g\|_{C^1}\leq C,
\ee
where the constant $C$ is independent of $\al$ and $\Phi$.
\end{lemma}
\begin{proof}
First, choose $g_0$ satisfying
\be\label{jh2}
\left\{
\begin{aligned}
&g_0^{\pr\pr\pr}+4\al^2g_0^\pr=0,\\
&g_0(\pm1)=0,\q \int^1_{-1}g_0(\th) d\th=1. \\
\end{aligned}
\right.
\ee
It has the following representation formula
\bn
g_0(\th)=\frac{\alpha(\cos (2 \alpha \th)-\cos (2 \alpha))}{\sin (2 \alpha)-2 \alpha \cos (2 \alpha)},
\en
and satisfies
\bn
\|g_0(\th)\|_{C^1}\ls 1.
\en
Set
\bn
h\ed g-g_0\,,
\en
Then $h$ satisfies
\be\label{jh3}
\left\{
\begin{aligned}
&h^{\pr\pr\pr}+4\al^2 h^\pr=-2\al\Phi(h h^{\prime}+h g_0^{\prime}+g_0 h^{\prime}+g_0 g_0^{\prime})\,,\\
&h(\pm1)=0,\q \int^1_{-1}h(\th) d\th=0\,. \\
\end{aligned}
\right.
\ee
We consider the corresponding linearized solution  $h_{Lin}$, which satisfies
\be\label{jh4}
\left\{
\begin{array}{l}
h_{Lin}^{\pr\pr\pr}+4\al^2 h^\pr_{Lin}=F(\th),\\
h_{Lin}(\pm1)=0,\q \int^1_{-1}h_{Lin}(\th) d\th=0. \\
\end{array}
\right.
\ee
The general solution of \eqref{jh4}$_1$ is given by
\be\label{jh5}
h_{Lin}(\th)= \mathcal{A}+ \mathcal{B}\cos(2\al\th)+\mathcal{C}\sin(2\al\th)+\f{1}{2\al}\int^\th_{-1}\int^{\bar{\th}}_{-1} F(s)\sin(2\al(\bar{\th}-s))dsd\bar{\th}.
\ee
From the boundary conditions and the integration condition, we have
\[
\left\{
\begin{aligned}
&\mathcal{A}+ \mathcal{B}\cos(2\al)+\mathcal{C}\sin(2\al)=-\f{1}{2\al}\int^1_{-1}\int^{\bar{\th}}_{-1} F(s)\sin(2\al(\bar{\th}-s))dsd\bar{\th},\\
&\mathcal{A}+ \mathcal{B}\cos(2\al)-\mathcal{C}\sin(2\al)=0,\nn\\
&2\mathcal{A}+\mathcal{B}\f{\sin(2\al)}{\al}=-\f{1}{2\al}\int^1_{-1}\int^\th_{-1}\int^{\bar{\th}}_{-1} F(s)\sin(2\al(\bar{\th}-s))dsd\bar{\th}d\th.
\end{aligned}
\right.
\]
Then we obtain
\[
\begin{aligned}
\mathcal{A}=&-\f{\sin(2\al)}{4\al(\sin(2\al)-2\al\cos(2\al))}\int^1_{-1}\int^{\bar{\th}}_{-1} F(s)\sin(2\al(\bar{\th}-s))dsd\bar{\th}\nn\\
   &+\f{\cos(2\al)}{2(\sin(2\al)-2\al\cos(2\al))}\int^1_{-1}\int^\th_{-1}\int^{\bar{\th}}_{-1} F(s)\sin(2\al(\bar{\th}-s))dsd\bar{\th}d\th,\\
\mathcal{B}=&\f{1}{2(\sin(2\al)-2\al\cos(2\al))}\lt\{ \int^1_{-1}\int^{\bar{\th}}_{-1} F(s)\sin(2\al(\bar{\th}-s))dsd\bar{\th}\rt.\\
            &\lt.\qq\qq\qq\qq\qq\q-\int^1_{-1}\int^\th_{-1}\int^{\bar{\th}}_{-1} F(s)\sin(2\al(\bar{\th}-s))dsd\bar{\th}d\th\rt\},\nn\\
\mathcal{C}=&-\f{1}{4\al\sin(2\al)}\int^1_{-1}\int^{\bar{\th}}_{-1} F(s)\sin(2\al(\bar{\th}-s))dsd\bar{\th}.
\end{aligned}
\]
Direct calculations imply that
\begin{align}
|\mathcal{A}|+|\mathcal{B}|\leq C\al^{-2}\|F\|_{L^\i},\q |\mathcal{C}|\leq C\al^{-1}\|F\|_{L^\i}, \label{jh9}
\end{align}
where the constant $C$ is independent of $\al$ and $\Phi$.
Direct calculations imply that
\begin{align}
&h_{Lin}^\pr(\th)=-2\al \mathcal{B}\sin(2\al\th)+2\al \mathcal{C}\cos(2\al\th)+\f{1}{2\al}\int^{{\th}}_{-1} F(s)\sin(2\al({\th}-s))ds, \label{jh6}\\
&h_{Lin}^{\pr\pr}(\th)=-(2\al)^2 \mathcal{B}\cos(2\al\th)-(2\al)^2 \mathcal{C}\sin(2\al\th)+\int^{{\th}}_{-1} F(s)\cos(2\al({\th}-s))ds,\label{jh7}\\
&h_{Lin}^{\pr\pr\pr}(\th)=(2\al)^3 \mathcal{B}\sin(2\al\th)-(2\al)^3 \mathcal{C}\cos(2\al\th)+ F(\th)-(2\al)\int^\th_{-1}F(s)\sin(2\al({\th}-s))ds.\label{jh8}
\end{align}

Now we use the contraction mapping principle to obtain a solution of \eqref{jh3}. We discuss the problem in two cases.

\textbf{Case I: $\al>\Phi^{1/3}$.}

Given $\t{h}\in C^1$, we define the map
\bn
T_>: \t{h}\longrightarrow h,
\en
where $h=T_>\t{h}$ solves the following boundary value problem
\[
\left\{
\begin{aligned}
&h^{\pr\pr\pr}+4\al^2 h^\pr=-2\al\Phi(\t{h}\t{h}^{\prime}+\t{h} g_0^{\prime}+g_0 \t{h}^{\prime}+g_0 g_0^{\prime})\,,\\
&h(\pm1)=0,\q \int^1_{-1}h(\th) d\th=0. \\
\end{aligned}\right.
\]
We first show that the map $T_>$ is well defined: From \eqref{jh5}, we have
\bn
h(\th)= \mathcal{A}+ \mathcal{B}\cos(2\al\th)+\mathcal{C}\sin(2\al\th)+\f{1}{2\al}\int^\th_{-1}\int^{\bar{\th}}_{-1} F_{\t{h}}(s)\sin(2\al(\bar{\th}-s))dsd\bar{\th},
\en
where
\bn
F_{\t{h}}(\th)=-2\al\Phi(\t{h}\t{h}^{\prime}+\t{h} g_0^{\prime}+g_0 \t{h}^{\prime}+g_0 g_0^{\prime})\,.
\en
Direct calculation shows
\be\label{jh10}
\begin{aligned}
|F_{\t{h}}(\th)|&\leq C\al\Phi\lt(|\t{h}||\t{h}^\pr|+|\t{h}|+|\t{h}^\pr|+1\rt)\\
&\leq C\al\Phi\lt(\|\t{h}\|_{L^\i}\|\t{h}^\pr\|_{L^\i}+\|\t{h}\|_{L^\i}+\|\t{h}^\pr\|_{L^\i}+1\rt)\,. \end{aligned}
\ee
Here the constant $C$ is independent of $\al$ and $\Phi$. From \eqref{jh5}, \eqref{jh6}, \eqref{jh7} and \eqref{jh8}, we see that
\[
\begin{aligned}
\|h\|_{L^\i}&\leq C\lt(|\mathcal{A}|+|\mathcal{B}|+|\mathcal{C}|\rt)+\|F_{\t{h}}\|_{L^\i},\\
\|h^{\pr}\|_{L^\i}&\leq C\lt(\al^2|\mathcal{B}|+\al|\mathcal{C}|\rt)+\|F_{\t{h}}\|_{L^\i}\,.
\end{aligned}
\]
where the constant $C$ is independent of $\al$ and $\Phi$.
From \eqref{jh9} and the above estimates we have
\bn
\|h\|_{L^\i}\leq C\al^{-2}\|F_{\t{h}}\|_{L^\i},\q \|h^\pr\|_{L^\i}\leq C\|F_{\t{h}}\|_{L^\i}.
\en
Moreover, in view of \eqref{jh10}, we deduce that
\be\label{jh1212}
\begin{aligned}
\|h\|_{L^\i}&\leq C\f{\Phi}{\al}\lt(\|\t{h}\|_{L^\i}\|\t{h}^\pr\|_{L^\i}+\|\t{h}\|_{L^\i}+\|\t{h}^\pr\|_{L^\i}+1\rt)\,;\\
\|h^\pr\|_{L^\i}&\leq C\al\Phi\lt(\|\t{h}\|_{L^\i}\|\t{h}^\pr\|_{L^\i}+\|\t{h}\|_{L^\i}+\|\t{h}^\pr\|_{L^\i}+1\rt)\,.\\
\end{aligned}
\ee
Then one concludes from \eqref{jh1212} that
\ba\l{jh12}
\|h\|_{C^1}\leq& C\al\Phi(1+(\al\Phi)^{-\f{1}{2}})\left(\|\t{h}\|_{C^1}^2+\|\t{h}\|_{C^1}+1\right)\\
\leq&C(\al\Phi)^{\f{1}{2}}\left(\|\t{h}\|_{C^1}^2+\|\t{h}\|_{C^1}+1\right)\,.
\ea
Here we assume $C>5$ without loss of generality. This shows that when $ 0\leq\al\Phi\leq \big(\f{1}{4C}\big)^2$ and $\|\t{h}\|_{C^1}\leq2C(\al\Phi)^{\f{1}{2}}$, one has
\bn
\|h\|_{C^1}\leq&C(\al\Phi)^{\f{1}{2}}\left(4C^2\al\Phi+2C(\al\Phi)^\f{1}{2}+1\right)<2C(\al\Phi)^{\f{1}{2}}\,.
\en
This shows that the mapping
\bn
T_>: \t{h}\longrightarrow h,
\en
maps the ball $B:=\{\|h\|_{C^1}\leq 2C(\al\Phi)^{\f{1}{2}}\}$ into itself. Also, if $\t{h}_1,\t{h}_2\in B$, by following the estimates in \eqref{jh12} line by line, we obtain
\bn
\|h_1-h_2\|_{C^1}\leq C(\al\Phi)^{\f{1}{2}}\|\t{h}_1-\t{h}_2\|_{C^1}\|(\t{h}_1,\t{h}_2)\|_{C^1}+C(\al\Phi)^{\f{1}{2}}\|\t{h}_1-\t{h}_2\|_{C^1}\leq\f{3}{8}\|\t{h}_1-\t{h}_2\|_E\,,
\en
which shows that $T$ is a contraction. Then we have a unique $h\in B$ such that
\bn
T_>h=h.
\en
In this case, we also have $\|h\|_{C^1}\ls 1$. This shows that the system \eqref{jh3} has a solution satisfying \eqref{jh2exg}, which also implies that the system \eqref{jh2} has a solution satisfying \eqref{jh2ex}.

\textbf{Case II: $\al\leq\Phi^\f{1}{3}$.}

In this case, we define the map
\bn
T_<: \t{h}\longrightarrow h,
\en
where $h=T_<\t{h}$ solves the following boundary value problem
\[
\left\{
\begin{aligned}
&h^{\pr\pr\pr}=-4\al^2 \t{h}^\pr-2\al\Phi(\t{h}\t{h}^{\prime}+\t{h} g_0^{\prime}+g_0 \t{h}^{\prime}+g_0 g_0^{\prime})\,,\\
&h(\pm1)=0,\q \int^1_{-1}h(\th) d\th=0. \\
\end{aligned}\right.
\]
Direct calculation shows
\ba\l{EH0415}
h(\th)&=a+b\th+c\th^2+\f{1}{2}\int_0^\th(\th-s)^2G_{\t{h}}(s)ds\,,\\
h^\pr(\th)&=b+2c\th+\int_0^\th(\th-s)G_{\t{h}}(s)ds\,.
\ea
where
\bn
a=&\frac12\int_0^1(1-s)^2G_{\t{h}}(s)ds-\frac12\int_{-1}^0(1+s)^2G_{\t{h}}(s)ds-\frac14\int_0^1(1-s)^3G_{\t{h}}(s)ds\\
&-\frac14\int_{-1}^0(1+s)^3G_{\t{h}}(s)ds\,;\\
b=&-\frac14\int_0^1(1-s)^2G_{\t{h}}(s)ds-\frac14\int_{-1}^0(1+s)^2G_{\t{h}}(s)ds\,;\\
c=&-\frac34\int_0^1(1-s)^2G_{\t{h}}(s)ds+\frac34\int_{-1}^0(1+s)^2G_{\t{h}}(s)ds+\frac14\int_0^1(1-s)^3G_{\t{h}}(s)ds\\
&+\frac14\int_{-1}^0(1+s)^3G_{\t{h}}(s)ds\,.
\en
Here
\[
G_{\t{h}}(s)\ed-4\al^2 \t{h}^\pr-2\al\Phi(\t{h}\t{h}^{\prime}+\t{h} g_0^{\prime}+g_0 \t{h}^{\prime}+g_0 g_0^{\prime})\,,
\]
and clearly
\ba\l{EG0415}
\|G_{\t{h}}\|_{L^\i}\les \al^2\|\t{h}^\pr\|_{L^\i}+\al\Phi\lt(\|\t{h}\|_{L^\i}\|\t{h}^\pr\|_{L^\i}+\|\t{h}\|_{L^\i}+\|\t{h}^\pr\|_{L^\i}+1\rt)\,.
\ea
Thus one concludes from \eqref{EH0415} and \eqref{EG0415} that
\bn
\|h\|_{C^1}\leq & C\al^2\|\t{h}^\pr\|_{L^\i}+C\al\Phi\lt(\|\t{h}\|_{L^\i}\|\t{h}^\pr\|_{L^\i}+\|\t{h}\|_{L^\i}+\|\t{h}^\pr\|_{L^\i}+1\rt)\\
\leq&C(\al^2+\al\Phi)\left(\|\t{h}\|_{C^1}^2+\|\t{h}\|_{C^1}+1\right)\,.
\en
The remaining proof in this case is the same as in Case I, and we therefore omit the details. In this case, we also have $\|h\|_{C^1}\ls \al^2+\al\Phi\ls 1$. This finishes the proof of the lemma.
\end{proof}

\section{$H^1$ weak formulation of system (\ref{nslog}) }\l{SEC3}

In this section, we study $H^1$ weak solutions of system \eqref{nslog}. For simplicity, we rewrite \eqref{nslog} as follows
\be\label{nslog2}
\left\{
\begin{aligned}
&e^{s}\lt[(u^{\e}\p_s+v^{\e}\p_\th) u^{\e} -\left(v^{\e}\right)^2\rt]+e^s\p_s p^{\e}- \left(\Dl_s u^\e-2\p_\th v^{\e}-{u^{\e} }\right)=0, \q \text{in}\q \O^\e,\\
&e^s\lt[(u^{\e}\p_s+v^{\e}\p_\th) v^{\e} +u^{\e} v^{\e}\rt]+e^s\p_\th p^{\e}-\left(\Dl_s v^\e+2 \p_\th u^{\e}-{v^{\e}}\right)=0,\hskip.2cm\q \text{in}\q \O^\e,\\
&\p_s (e^su^{\e})+\p_\th (e^sv^\e)=0,  \q \text{in}\q \O^\e,\\
&\lt(u^{\e},v^\e\rt)=0, \q \text{on}\q \p\O^\e,\\
&\int^{\al+\th^\e_+(s)}_{-\al-\th^\e_{-}(s)} e^su^\e(\th) d\th=\Phi,\q \text{ for  a.e. } s\in\bR.
\end{aligned}
\right.
\ee
We now define weak solutions of \eqref{nslog2}.

\begin{definition}\l{Def1}
We call
\[
(u^\e,v^\e)\in \mathbf{H}^1_{0,\s}(\O^\e)\ed\lt\{(u,v)\in \lt(H^1_0(\O^\e)\rt)^2\,:\,\p_s(e^su)+\p_\th(e^sv)=0\rt\}
\]
a weak solution of \eqref{nslog2} if the following conditions hold:
\begin{itemize}
  \item[(i)] for any pair of test functions
\[
(\phi,\psi)\in \mathbf{C}^\i_{c,\s}(\O^\e)\ed\lt\{(\phi,\psi)\in \lt(C_c^\infty(\O^\e)\rt)^2\,:\,\p_s(e^s\phi)+\p_\th(e^s\psi)=0\rt\}\,,
\]
we have the following identity:
\begin{align}\l{Weak1}
&\int_{\O^\e} \na u^\e\cdot \na \phi + \na v^\e\cdot \na \psi- (-2\p_\th v^\e-u^\e)\phi- (2\p_\th u^\e-v^\e)\psi dx \nn\\
&+\int_{\O^\e}e^{s}\lt[(u^{\e}\p_s+v^{\e}\p_\th) u^{\e} -\left(v^{\e}\right)^2\rt]\phi dx\nn\\
&+ \int_{\O^\e}e^s\lt[(u^{\e}\p_s+v^{\e}\p_\th) v^{\e} +u^{\e} v^{\e}\rt]\psi dx=0\,;
\end{align}
  \item[(ii)] $u^\e$ satisfies the flux condition \eqref{nslog2}$_5$.
\end{itemize}
\end{definition}

In this paper, we focus on solutions of \eqref{nslog2} around the Jeffery--Hamel flow
\bn
(u^0,v^0)=(e^{-s} \t{f}(\th), 0), \q (s,\th)\in \O^\e\,,
\en
with
\[
\t{f}(\th)=\left\{
\begin{aligned}
&f(\th)\,,&\q\text{for}\q\th\in(-\al,\al)\,;\\
&0\,,&\text{else}\,.
\end{aligned}
\right.
\]
Here the existence of $f(\th)$ is given by Proposition \ref{jh1prop}. Therefore it is more natural to consider the difference between $u^\e$ and the underlying background solution:
\bn
(u,v)\ed(u^\e,v^\e)-(u^0,v^0)\,.
\en
Subtracting the Jeffery--Hamel flow from \eqref{nslog2}, we see that $(u,v)$ satisfies
\be\label{nslogerror}
\left\{
\begin{aligned}
&e^{s}\lt[(u\p_s+v\p_\th) u+ (u\p_s+v\p_\th) u^0+ (u^0\p_s+v^0\p_\th) u-v^2\rt]+e^s\p_s p\\
&\hskip.5cm-\left(\Dl_s u-2\p_\th v-u\right)=0\,,\q\text{in}\q \O^\e/\Sigma_\pm\,,\\
&e^s\lt[(u\p_s+v\p_\th) v+ (u\p_s+v\p_\th) v^0+ (u^0\p_s+v^0\p_\th) v +(u+u^0) v\rt]+e^s\p_\th p\\
&\hskip.5cm-\left(\Dl_s v+2 \p_\th u-{v}\right)=0\,,\q\text{in}\q \O^\e/\Sigma_\pm\,,\\
&\p_s (e^su)+\p_\th (e^sv)=0\,,\q\text{in}\q \O^\e/\Sigma_\pm,\\
&\lt(u,v\rt)=0\,, \q \text{on}\q \p\O^\e\,,\\
&\int^{\al+\th^\e_+(s)}_{-\al-\th^\e_{-}(s)} e^su(\th) d\th=0\,,\q \text{for  a.e. } s\in\bR\,.
\end{aligned}
\right.
\ee
Here $p= p^\e-p^0$, where
\[
p^{0}\ed\lt\{
\bali
& 2 e^{-2s}f(\th)-\f{1}{2}e^{-2s} C_{\Phi,\al}, \q (s,\th)\in \O;\\
& 2 e^{-2s}f(\pm\al)-\f{1}{2}e^{-2s} C_{\Phi,\al},\q (s,\th)\in \O^\e_\pm,
\eali
\rt.
\]
is the natural continuous extension of \eqref{jefferyhamelsol1}$_3$ in $\o^\e$. Meanwhile, the following compatibility condition on $\Sigma_{\pm}\ed\{\th=\pm \al, s\in \bR \}$ holds
\be\label{compatibility1}
\left\{
\begin{aligned}
&[\p_\th  u]\big|_{\Sigma_\pm}=\pm e^{-s} f'(\pm \al),\\
&[\p_\th v-e^sp]\big|_{\Sigma_\pm}=0.
\end{aligned}
\right.
\ee
Here $[g]\big|_{\Sigma_\pm}\ed g(s,\pm\al_+)-g(s,\pm\al_-)$ is the jump of a function $g$ on $\O^\e$ crossing lines $\th=\pm\al$, respectively.

The following is the weak formulation of \eqref{nslogerror}--\eqref{compatibility1}.
\begin{definition}\l{Def2}
 We call $(u,v)\in \mathbf{H}^1_{0,\s}(\O^\e)$ a weak solution of \eqref{nslogerror}--\eqref{compatibility1}, if and only if
\begin{itemize}
  \item[(i)] for any pair of test functions $(\phi,\psi)\in \mathbf{C}^\i_{c,\s}(\O^\e)$, we have the following identity:
\bn
&\int_{\O^\e} \p_su\p_s\phi+ \p_sv\p_s\psi+(\p_\th v+u)(\p_\th\psi+\phi)+(\p_\th u-v)(\p_\th \phi-\psi)\,d\th d s\\
+&f^\pr(\al)\int_{-\i}^\i e^{-s}\phi(s,\al)ds-f^\pr(-\al)\int_{-\i}^\i e^{-s}\phi(s,-\al)ds\\
+&\int_{\O^\e}e^s\phi\lt((u\p_s+v\p_\th)u-v^2\rt)+e^s\psi\lt((u\p_s+v\p_\th)v+uv\rt)d\th ds\\
+& \int_{\O^\e}\t{f}(\th)\lt((\p_su-u)\phi+(\p_sv+v)\psi\rt)+\t{f}^\pr(\th)v\phi\,d\th ds =0\,;
\en
  \item[(ii)] $u$ satisfies the following flux condition
  \be\l{VFlux}
\int^{\al+\th^\e_+(s)}_{-\al-\th^\e_{-}(s)} e^su(s,\th) d\th=0\,,\q \text{ for  any } s\in\bR\,.
\ee
\end{itemize}

\end{definition}

Now we give a brief explanation of the condition \eqref{compatibility1}. Suppose $(u^\e,v^\e)$ is smooth enough. Multiplying \eqref{nslog2}$_{1,2}$ by $(\phi,\psi)\in \mathbf{C}^\i_{c,\s}(\O^\e)$ and integrating the resulting equations on $\O^\e$, we see that
\begin{align}
&\int_{\O^\e}- \Dl u^\e \phi - \Dl v^\e \psi- (-2\p_\th v^\e-u^\e)\phi- (2\p_\th u^\e-v^\e)\psi dx+\int_{\O^\e} e^s\p_s p^\e \phi+e^s\p_\th p^\e \psi dx\nn\\
&+\int_{\O^\e}e^{s}\lt[(u^{\e}\p_s+v^{\e}\p_\th) u^{\e} -\left(v^{\e}\right)^2\rt]\phi dx+ \int_{\O^\e}e^s\lt[(u^{\e}\p_s+v^{\e}\p_\th) v^{\e} +u^{\e} v^{\e}\rt]\psi dx =0.\nn
\end{align}
Integration by parts indicates that
\be\l{042201}
\begin{aligned}
&\int_{\O^\e} \na u^\e\cdot \na \phi +\na v^\e\cdot\na \psi- (-2\p_\th v^\e-u^\e)\phi- (2\p_\th u^\e-v^\e)\psi dx\\
&+\int_{\Sigma_\pm} [\p_\th u^\e] \phi ds+[\p_\th v^\e-e^s p^\e] \psi ds\\
&+\int_{\O^\e}e^{s}\lt[(u^{\e}\p_s+v^{\e}\p_\th) u^{\e} -\left(v^{\e}\right)^2\rt]\phi dx\\
&+ \int_{\O^\e}e^s\lt[(u^{\e}\p_s+v^{\e}\p_\th) v^{\e} +u^{\e} v^{\e}\rt]\psi dx=0.
\end{aligned}
\ee
Recall \eqref{Weak1} in Definition \ref{Def1}, the second line in \eqref{042201} must vanish. Thus we need
\be\label{compatibility}
[\p_\th u^\e]\big|_{\Sigma_\pm}=0, \q [\p_\th v^\e-e^s p^\e]\big|_{\Sigma_\pm}=0\,.
\ee
This gives \eqref{compatibility1}.

\section{Solvability of system (\ref{nslogerror}) with the compatibility condition (\ref{compatibility1})}\l{SEC4}
In this section, we prove the existence of an $H^1$ weak solution of system \eqref{nslog}. We first restate Theorem \ref{thmexistence1} by using Definition \ref{Def2}.
\begin{theorem}[Restatement of Theorem \ref{thmexistence1}]\label{thmexistence}
There exists a constant $C_\ast>0$, independent of $\al$, $\Phi$, $\kappa$, and $\e$, such that when $\e\leq\min\{\f{\pi-2\al}{4}, \al, \f{\kappa}{4}\}$, and $\al$, $\Phi$ satisfy
\be\label{condition0}
\al\Phi + \f{\Phi}{(\pi-2\al)^2}\leq C_\ast\,,
\ee
the system \eqref{nslogerror} with the compatibility condition \eqref{compatibility1} has a weak solution satisfying
\ba\label{UBDxx}
&\|\p_s(u,v)\|_{L^2(\O^\e)}^2+\f{\lt(\pi-2\al\rt)^2}{32\pi^2}\|\p_\th (u,v)\|_{L^2(\O^\e)}^2\leq C\f{\Phi^2\e^2}{\al^4\kappa(\pi-2\al)^2}\,.
\ea
\end{theorem}

The proof proceeds in three steps. We first reformulate the original problem on a truncated bounded domain $\O^\e_R$, and then establish the solvability of the truncated problem by applying the Leray–Schauder fixed-point theorem. Finally, by letting the truncated domains exhaust $\O^\e$ and passing to the limit, we complete the proof of Theorem \ref{thmexistence}.

\subsection{Weak solutions in a truncated domain}

Our proof in this section is based on the Leray--Schauder fixed point theorem. In order to deal with the lack of compactness of the Sobolev embedding on unbounded domains, we first study the corresponding problem on the truncated domain
\bn
\O^\e_R\ed\lt\{(s,\th)\in\O^\e\,:\,-R<s<R\rt\}\,.
\en
The corresponding solution space is
\bn
\mathbf{H}^1_{0,\s}(\O^\e_R)\ed\lt\{(u_R,v_R)\in \lt(H^1_0(\O^\e_R)\rt)^2\,:\,\p_s(e^su_R)+\p_\th(e^sv_R)=0\rt\}\,.
\en
Since, by the zero boundary condition on $s=R$ and the divergence-free property of $(u_R,v_R)$, the vanishing flux condition \eqref{VFlux} holds automatically. Therefore, according to Definition \ref{Def2},  the existence of a weak solution in the truncated domain $\O^\e_R$ can be described as follows:
\begin{problem}\l{P1}
 Find $(u_R,v_R)\in \mathbf{H}^1_{0,\s}(\O^\e_R)$ such that for any test function pair
\[
(\phi,\psi)\in \mathbf{C}^\i_{c,\s}(\O^\e_R)\ed\lt\{(\phi,\psi)\in \lt(C_c^\infty(\O^\e_R)\rt)^2\,:\,\p_s(e^s\phi)+\p_\th(e^s\psi)=0\rt\}\,,
\]
the following identity holds:
 \ba\l{EP1}
&\int_{\O^\e_R} \p_su_R\p_s\phi+ \p_sv_R\p_s\psi+(\p_\th v_R+u_R)(\p_\th\psi+\phi)+(\p_\th u_R-v_R)(\p_\th \phi-\psi)\,d\th d s\\
+&f^\pr(\al)\int_{-R}^R e^{-s}\phi(s,\al)ds-f^\pr(-\al)\int_{-R}^R e^{-s}\phi(s,-\al)ds\\
+&\int_{\O^\e_R}e^s\phi\lt((u_R\p_s+v_R\p_\th)u_R-v_R^2\rt)+e^s\psi\lt((u_R\p_s+v_R\p_\th)v_R+u_Rv_R\rt)d\th ds\\
+& \int_{\O^\e_R}\t{f}(\th)\lt((\p_su_R-u_R)\phi+(\p_sv_R+v_R)\psi\rt)+\t{f}^\pr(\th)v_R\phi\,d\th ds =0\,.
\ea
\end{problem}

\qed

Before proving solvability, we state two lemmas that will be used in the subsequent estimates.
\begin{lemma}[Wirtinger's inequality, in Section 7.7 of \cite{Hardy1948}]\label{lemwir}
Let $g$ be a continuously differentiable function on $[-\beta,\beta]$ with $g(-\beta)=g(\beta)=0$. Then
\[
\int_{-\beta}^{\beta} g^2(\theta)d\theta\leq\f{4\beta^2}{\pi^2}\int_{-\beta}^{\beta}|g'(\theta)|^2d\theta\,.
\]
The constant $4\beta^2/\pi^2$ is optimal, and the equality holds if and only if $g(\theta)=C\cos(\f{\pi}{2\beta}\theta)$, where $C$ is an arbitrary constant.
\end{lemma}

\begin{lemma}[Compound Poincar\'e inequality]\l{LemCP}
Let $(u,v)\in \lt(H^1_0(\O^\e_R)\rt)^2$, then
\bn
\|u\|_{L^2(\O^\e_R)}+\|v\|_{L^2(\O^\e_R)}\leq\f{2\al+2\ve}{\pi-(2\alpha+2\ve)}\left(\|\p_\theta v+u\|_{L^2(\O^\e_R)}+\|\p_\theta u-v\|_{L^2(\O^\e_R)}\right)\,.
\en
\end{lemma}

\begin{proof}
Define $\mu\ed\f{\al+2\e}{\pi}$, using Lemma \ref{lemwir}, we find
\be\label{CP1}
\begin{split}
\left(\int_{\O^\e_R}u^2dsd\th\right)^{1/2}&\leq \mu\left(\int_{\O^\e_R}(\p_\theta u)^2dsd\th\right)^{1/2}\\
&\leq\mu\left(\int_{\O^\e_R}(\p_\theta u-v)^2dsd\th\right)^{1/2}+\mu\left(\int_{\O^\e_R}v^2dsd\th\right)^{1/2}\,.
\end{split}
\ee
Meanwhile,
\be\label{CP2}
\begin{split}
\left(\int_{\O^\e_R}v^2dsd\th\right)^{1/2}&\leq \mu\left(\int_{\O^\e_R}(\p_\theta u)^2dsd\th\right)^{1/2}\\
&\leq\mu\left(\int_{\O^\e_R}(\p_\theta v+u)^2 dsd\th\right)^{1/2}+\mu\left(\int_{\O^\e_R}u^2 dsd\th\right)^{1/2}\,.
\end{split}
\ee
Combining \eqref{CP1} and \eqref{CP2}, one derives
\[
\|u\|_{L^2(\O^\e_R)}+\|v\|_{L^2(\O^\e_R)}\leq\f{\mu}{1-\mu}\left(\|\p_\theta v+u\|_{L^2(\O^\e_R)}+\|\p_\theta u-v\|_{L^2(\O^\e_R)}\right)\,.
\]
This concludes the proof of the lemma.
\end{proof}

\subsection{Solvability of the truncated problem}
For convenience, we define the linear operator $\mL_R\,:\,\mathbf{H}^1_{0,\s}(\O^\e_R)\,\to\,\lt(\mathbf{H}^1_{0,\s}(\O^\e_R)\rt)^\pr$ that
\ba\l{DLR}
&\langle \mL_R(u_R,v_R)\,,\,(\phi,\psi)\rangle\\
\ed&\int_{\O^\e_R} \p_su_R\p_s\phi+ \p_sv_R\p_s\psi+(\p_\th v_R+u_R)(\p_\th\psi+\phi)+(\p_\th u_R-v_R)(\p_\th \phi-\psi)\,d\th d s\,.
\ea
The following proposition gives the key property of the operator $\mL_R$.
\begin{proposition}\l{Propinv}
$\mL_R$ is continuous and invertible.
\end{proposition}
\begin{proof}
Clearly $\mL_R$ is well defined and continuous. To prove $\mL_R$ is invertible, it remains to show $\langle\mL_R\cd\,,\,\cd\rangle$ is coercive. Indeed, by \eqref{DLR}
\bn
\langle \mL_R(u_R,v_R)\,,\,(u_R,v_R)\rangle=&\int_{\O^\e_R} \lt( (\p_s u_R)^2+(\p_s v_R)^2+ (\p_\th u_R-v_R)^2+(\p_\th v_R+u_R)^2 \rt)\,d\th ds\,.
\en
Using Minkowski's inequality, we have
\ba\l{CPI1}
&\|\p_\th u_R\|_{L^2(\O^\e_R)}+\|\p_\th v_R\|_{L^2(\O^\e_R)}\\
\leq&\|\p_\th u_R-v_R\|_{L^2(\O^\e_R)}+\|\p_\th v_R+u_R\|_{L^2(\O^\e_R)}+\|v_R\|_{L^2(\O^\e_R)}+\|u_R\|_{L^2(\O^\e_R)}\\
\leq&\f{\pi}{\pi-(2\al+2\ve)}\left(\|\p_\th u_R-v_R\|_{L^2(\O^\e_R)}+\|\p_\th v_R+u_R\|_{L^2(\O^\e_R)}\right)\,.
\ea
In the last line, we have applied Lemma \ref{LemCP}. This indicates
\ba\l{ECO}
&\int_{\O^\e_R}\lt((\p_\th u_R-v_R)^2+(\p_\th v_R+u_R)^2 \rt)\,d\th ds\\
&\geq \f{1}{2}\left(\|\p_\th u_R-v_R\|_{L^2(\O^\e_R)}+\|\p_\th v_R+u_R\|_{L^2(\O^\e_R)}\right)^2\\
&\geq \f{\lt(\pi-(2\al+2\ve)\rt)^2}{2\pi^2}\left(\|\p_\th u_R\|_{L^2(\O^\e_R)}^2+\|\p_\th v_R\|_{L^2(\O^\e_R)}^2\right)\,.
\ea
This completes the proof of this proposition.
\end{proof}
The nonlinear operator $\mN_R\,:\,\mathbf{H}^1_{0,\s}(\O^\e_R)\,\to\,\lt(\mathbf{H}^1_{0,\s}(\O^\e_R)\rt)^\pr$ is defined by
\ba\l{ENR}
&\langle \mN_R(u_R,v_R)\,,\,(\phi,\psi)\rangle\\
\ed&\int_{\O^\e_R}e^s\phi\lt((u_R\p_s+v_R\p_\th)u_R-v_R^2\rt)+e^s\psi\lt((u_R\p_s+v_R\p_\th)v_R+u_Rv_R\rt)d\th ds\\
&+\int_{\O^\e_R}\t{f}(\th)\lt((\p_su_R-u_R)\phi+(\p_sv_R+v_R)\psi\rt)+\t{f}^\pr(\th)v_R\phi\,d\th ds\,.
\ea
Clearly $\mN_R$ is well-defined and continuous. Finally, we denote $\mf_R\in\lt(\mathbf{H}^1_{0,\s}(\O^\e_R)\rt)^\pr$ that
\ba\l{EFFR}
\langle\mf_R\,,\,(\phi,\psi)\rangle\ed -f^\pr(\al)\int_{-R}^R e^{-s}\phi(s,\al)ds+f^\pr(-\al)\int_{-R}^R e^{-s}\phi(s,-\al)ds\,.
\ea
We \textbf{claim} this $\mf_R$ is well defined and its norm is independent of $R>1$. Indeed, using Proposition \ref{jh1prop} and Newton--Leibniz formula
\ba\l{FR1}
|\langle\mf_R\,,\,(\phi,\psi)\rangle|\leq& C \al^{-2}\Phi\lt|\int^{+R}_{-R} e^{-s}\lt|\int^{\pm\al\pm\th^\e_{\pm}(s)}_{\pm\al}\p_\th \phi(s,\th)d\th\rt| ds\rt|\\
     \leq&C\al^{-2}\Phi\lt|\int^{+R}_{-R} e^{-s}\lt(\int^{\pm\al\pm\th^\e_{\pm}(s)}_{\pm\al}(\p_\th \phi(s,\th))^2d\th\rt)^{1/2}|\th^\e_{\pm}(s)|^{1/2} ds\rt|\\
\leq&C\al^{-2}\Phi\int^{+R}_{-R} e^{-s}\lt(\int^{\pm\al\pm\th^\e_{\pm}(s)}_{\pm\al}(\p_\th \phi(s,\th))^2d\th\rt)^{1/2}\e^{1/2} e^{-\f{|s|\kappa}{2\e}} ds\,.
\ea
In the last line, we have applied the assumption of $\th^\e$ given in \eqref{Ass}. By choosing {$\e\leq\kappa/4$}, one has
\ba\l{FR2}
|\langle\mf_R\,,\,(\phi,\psi)\rangle|
     \leq&C { \al^{-2}\Phi}\sqrt{\e}\int^{+\i}_{-\i} \lt(\int^{\pm\al\pm\th^\e_{\pm}(s)}_{\pm\al}(\p_\th \phi(s,\th))^2d\th\rt)^{1/2} e^{-{ \f{\kappa}{4\e}}|s|} ds\\
     \leq& C{\al^{-2}\Phi}\sqrt{\e}\left(\int_{-\i}^\i e^{-\f{\kappa}{2\e}|s|}ds\right)^{1/2}\|\p_\th \phi\|_{L^2(\O^\e_R)}\\
\leq& C {\al^{-2}\Phi \kappa^{-1/2}\e }\|(\phi,\psi)\|_{\mathbf{H}^1_{0,\sigma}(\O^\e_R)}\,,
\ea
which proves the \textbf{claim}. Problem \ref{P1} is equivalent to $\mL_R(u_R,v_R)+\mN_R(u_R,v_R)=\mf_R$, which is $(u_R,v_R)=\mL^{-1}_R\lt(\mf_R-\mN_R(u_R,v_R)\rt)\ed \mT_R(u_R,v_R)$. Using the \emph{Leray-Schauder fixed point theorem}, the existence of a weak solution to Problem \ref{P1} follows from the following lemma:
\begin{lemma}
The operator $\mT_R$ satisfies
\begin{itemize}
  \item[(i)] $\mT_R$ is compact;
  \item[(ii)] Any solution of $(u_R,v_R)=\la \mT_R(u_R,v_R)$ with $\la\in(0,1)$ have a uniform bound.
\end{itemize}
\end{lemma}
\begin{proof}
Let $\{(u_{R,k},v_{R,k})\}_{k=1}^\i\subset\mathbf{H}^1_{0,\s}(\O^\e_R)$ be a bounded sequence that
\bn
(u_{R,k},v_{R,k})\to(u_R,v_R)\q\text{weakly in}\q H^1(\O^\e_R)\,.
\en
 By the compactness of the Sobolev embedding, there exists a subsequence of $\{(u_{R,k},v_{R,k})\}$ (still denoted by $\{(u_{R,k},v_{R,k})\}$) that
\bn
(u_{R,k},v_{R,k})\to(u_R,v_R)\q\text{strongly in}\q L^2(\O^\e_R)\cap L^4(\O^\e_R)\,.
\en
In view of \eqref{ENR}, direct calculation shows
\ba\l{EN1}
&\langle\mN_R(u_{R,k},v_{R,k})-\mN_R(u_R,v_R)\,,\,(\phi,\psi)\rangle\\
=&\int_{\O^\e_{R,k}}e^s\phi\lt((u_{R,k}\p_s+v_{R,k}\p_\th)u_{R,k}-v_{R,k}^2\rt)+e^s\psi\lt((u_{R,k}\p_s+v_{R,k}\p_\th)v_{R,k}+u_{R,k}v_{R,k}\rt)d\th ds\\
&+\int_{\O^\e_R}\t{f}(\th)\lt(-u_{R,k}(\p_s\phi+\phi)+v_{R,k}(\psi-\p_s\psi)\rt)+\t{f}^\pr(\th)v_{R,k}\phi\,d\th ds\\
&-\int_{\O^\e_R}e^s\phi\lt((u_R\p_s+v_R\p_\th)u_R-v_R^2\rt)+e^s\psi\lt((u_R\p_s+v_R\p_\th)v_R+u_Rv_R\rt)d\th ds\\
&-\int_{\O^\e_R}\t{f}(\th)\lt(-u_R(\p_s\phi+\phi)+v_R(\psi-\p_s\psi)\rt)+\t{f}^\pr(\th)v_R\phi\,d\th ds\,.
\ea
A standard argument, based on adding and subtracting the same auxiliary term, together with H\"older's inequality, gives
\ba\l{EN2}
&|\langle\mN_R(u_{R,k},v_{R,k})-\mN_R(u_R,v_R)\,,\,(\phi,\psi)\rangle|\\[2mm]
\ls& e^R\|(u_{R,k},v_{R,k})-(u_R,v_R)\|_{L^4}\|(u_{R,k},v_{R,k})\|_{H^1}\|(\phi,\psi)\|_{H^1}\\[2mm]
&\hskip 1cm+e^R\|(u_{R,k},v_{R,k})-(u_R,v_R)\|_{L^4}\|(u_{R},v_{R})\|_{H^1}\|(\phi,\psi)\|_{H^1}\\[2mm]
&\hskip 1cm+{\al^{-2}\Phi}\|(u_{R,k},v_{R,k})-(u_R,v_R)\|_{L^2}\|(\phi,\psi)\|_{H^1}\\[2mm]
&\to 0\q\text{as}\q k\to\infty\,.
\ea
This implies
\bn
\mN_R(u_{R,k},v_{R,k})\to\mN_R(u_R,v_R)\q\text{strongly in}\q\lt(\mathbf{H}^1_{0,\s}(\O^\e_R)\rt)^\pr\,.
\en
In view of Proposition \ref{Propinv}, \[
\mL_R^{-1}\,:\,\lt(\mathbf{H}^1_{0,\s}(\O^\e_R)\rt)^\pr\to\mathbf{H}^1_{0,\s}(\O^\e_R)
\]
is continuous. This implies
\[
\mT_R(u_{R,k},v_{R,k})=\mL^{-1}_R\lt(\mf_R-\mN_R(u_{R,k},v_{R,k})\rt)\to\mL^{-1}_R\lt(\mf_R-\mN_R(u_{R},v_{R})\rt)=\mT_R(u_{R},v_{R})
\]
strongly in $\mathbf{H}^1_{0,\s}(\O^\e_R)$ as $k\to\i$. This proves part (i) of the lemma.

Moreover, equality $(u_R,v_R)=\la \mT_R(u_R,v_R)$ is equivalent to
\be\l{Homeq}
\mL_R(u_R,v_R)+\la \mN_R(u_R,v_R)=\la \mf_R\,.
\ee
Testing \eqref{Homeq} by $(u_R,v_R)$, recalling \eqref{EP1}, one has
\ba\l{EM0}
&\un{\int_{\O^\e_R} \lt( (\p_s u_R)^2+(\p_s v_R)^2+ (\p_\th u_R-v_R)^2+(\p_\th v_R+u_R)^2 \rt)\,d\th ds}_{M_1}\\
+&\un{\la\int_{\O^\e_R}\t{f}(\th)(v^2_R-u^2_R)+\t{f}^\pr(\th)v_Ru_R\,d\th ds}_{M_2}\\
+&\un{\la f^\pr(\al)\int_{-R}^R e^{-s}u_R(s,\al)ds-\la f^\pr(-\al)\int_{-R}^R e^{-s}u_R(s,-\al)ds}_{M_3}=0.
\ea
Using \eqref{ECO},
\ba\l{EM1}
M_1\geq \|\p_su_R\|_{L^2(\O^\e_R)}^2+\|\p_sv_R\|_{L^2(\O^\e_R)}^2+\f{\lt(\pi-(2\al+2\ve)\rt)^2}{2\pi^2}\left(\|\p_\th u_R\|_{L^2(\O^\e_R)}^2+\|\p_\th v_R\|_{L^2(\O^\e_R)}^2\right)\,.
\ea
Moreover, using integration by parts, \eqref{jh2ex} and Lemma \ref{lemwir}
\ba\l{EM2}
|M_2|\leq &C{\al^{-2}\Phi}\lt(\|u_R\|_{L^2(\O^\e_R)}^2+\|v_R\|_{L^2(\O^\e_R)}^2\rt)\\
\leq& C{\al^{-2}\Phi }\f{4(\al+\ve)^2}{\pi^2}\left(\|\p_\theta u_R\|_{L^2(\O^\e_R)}^2+\|\p_\theta v_R\|_{L^2(\O^\e_R)}^2\right)\,.
\ea
Finally, in view of \eqref{FR1}--\eqref{FR2} and using \eqref{CPI1}, one deduces
\begin{align}\l{EM3}
|M_3|\leq& C{\al^{-2}\Phi \kappa^{-1/2}\e}\|\p_\th u_R\|_{L^2(\O^\e_R)}\nn\\
\leq&\f{\lt(\pi-(2\al+2\ve)\rt)^2}{4\pi^2}\|\p_\th u_R\|_{L^2(\O^\e_R)}^2+C{\Phi^2}\f{\e^2}{\al^4\kappa(\pi-2(\al+\e))^2}\,.
\end{align}
Thus if
\be\label{conditionx}
C{\al^{-2}\Phi }\f{4(\al+\ve)^2}{\pi^2}<\f{\lt(\pi-(2\al+2\ve)\rt)^2}{8\pi^2}\,,
\ee
one obtains by substituting \eqref{EM1}--\eqref{EM2}--\eqref{EM3} in \eqref{EM0} that
\begin{align}\l{UBD}
&\|\p_su_R\|_{L^2(\O^\e_R)}^2+\|\p_sv_R\|_{L^2(\O^\e_R)}^2+\f{\lt(\pi-(2\al+2\ve)\rt)^2}{8\pi^2}\left(\|\p_\th u_R\|_{L^2(\O^\e_R)}^2+\|\p_\th v_R\|_{L^2(\O^\e_R)}^2\right)\nn\\
\leq& C{\f{\Phi^2\e^2}{\al^4\kappa(\pi-2(\al+\e))^2}}\,.
\end{align}
This holds uniformly for any $0<\la<1$. Here, since $\varepsilon\le\min\left\{\frac{\pi-2\alpha}{4},\,\alpha,\,\frac{\kappa}{4}\right\}$, we have
\[
\alpha+\varepsilon\le 2\alpha,
\qquad
\pi-2(\alpha+\varepsilon)\ge \frac{\pi-2\alpha}{2}.
\]
Therefore, the left-hand side of \eqref{conditionx} satisfies
\[
C\alpha^{-2}\Phi\frac{4(\alpha+\varepsilon)^2}{\pi^2}
\le
C\alpha^{-2}\Phi\frac{16\alpha^2}{\pi^2}
=
\frac{16C}{\pi^2}\Phi\,,
\]
while the right-hand side of \eqref{conditionx} admits
\[
\frac{(\pi-(2\alpha+2\varepsilon))^2}{8\pi^2}
=
\frac{(\pi-2(\alpha+\varepsilon))^2}{8\pi^2}
\ge
\frac{(\pi-2\alpha)^2}{32\pi^2}.
\]
Hence, a sufficient condition for \eqref{conditionx} is
\[
512C\Phi<(\pi-2\alpha)^2\,.
\]
Thus there exists a constant \(C_{*1}>0\), independent of \(\alpha\) and \(\Phi\), such that \eqref{conditionx} is ensured by
\ba\l{Cond1}
\frac{\Phi}{(\pi-2\alpha)^2}\le C_{*1}\,.
\ea
Moreover, by Proposition \ref{jh1prop}, the background Jeffery--Hamel flow exists provided
\ba\l{JFCond}
\al\Phi<\frak{c}_0\,.
\ea
Thus both \eqref{Cond1} and \eqref{JFCond} are satisfied provided
\[
\al\Phi+\frac{\Phi}{(\pi-2\alpha)^2}\le \min\{C_{*1},\frak{c}_0\}\,,
\]
which is exactly \eqref{condition0} by choosing $C_*=\min\{C_{*1},\frak{c}_0\}$. Under this condition, we have for some constant $C>0$, independent of $\al$, $\Phi$, $\kappa$ and $\e$,
\bn
\|\p_su_R\|_{L^2(\O^\e_R)}^2+\|\p_sv_R\|_{L^2(\O^\e_R)}^2+\f{\lt(\pi-2\al\rt)^2}{32\pi^2}\left(\|\p_\th u_R\|_{L^2(\O^\e_R)}^2+\|\p_\th v_R\|_{L^2(\O^\e_R)}^2\right)\leq& {\f{C\Phi^2\e^2}{\al^4\kappa(\pi-2\al)^2}}\,.
\en
This holds uniformly for any $0<\la<1$. This proves part (ii) and hence the lemma.
\end{proof}

\subsection{Solvability on the unbounded domain $\O^\e$}
Now it remains to pass to the limit as \(R \to \infty\) to finish the proof. In view of \eqref{UBD}, $(u_R,v_R)$ is uniformly bounded in $H^1(\O^\e)$. Therefore, by passing to a subsequence, we may assume that
\bn
&(u_R,v_R)\to(u,v)\,\q\text{weakly in}\q H^1(\O^\e)\q\text{as}\q R\to\i\,;\\
&(u_R,v_R)\to(u,v)\,\q\text{strongly in}\q L^4(\O^\e_{R_0})\q\text{as}\q R\to\i\q\text{for any }R_0>0\,.
\en
Thus for any $(\phi,\psi)\in \mathbf{C}_{c,\s}^\infty(\O^\e)$,
\[
\langle \mL_\infty(u_R,v_R)\,,\,(\phi,\psi)\rangle\to\langle \mL_\infty(u,v)\,,\,(\phi,\psi)\rangle\q\text{as}\q R\to\i\,.
\]
Meanwhile, there exists $R_0>0$ such that $\mathrm{supp}(\phi,\psi)\subset\O^\e_{R_0}$. Arguing as in \eqref{EN1}--\eqref{EN2}, we obtain
\[
\langle\mN_\infty(u_{R},v_{R})-\mN_\infty(u,v)\,,\,(\phi,\psi)\rangle\to 0\q\text{as}\q R\to\i\,.
\]
Recalling \eqref{EFFR}, we see that
\[
\langle\mf_R\,,\,(\phi,\psi)\rangle\to\langle\mf_\infty\,,\,(\phi,\psi)\rangle
\]
follows from the dominated convergence theorem. This completes the proof of the existence of a solution to \eqref{nslogerror}--\eqref{compatibility1}.

\section{More precise kinetic energy estimates in the middle strip $\O$}\l{SEC5}

For fixed $\al\in(0,\f{\pi}{2})$, $\Phi>0,\kappa>0$, from the energy estimate in \eqref{UBDxx} and Poincar\'e inequality, we see that the kinetic energy for the error velocity in the middle strip $\O$ is of order $\e^2$.  Next we will show that when the angle $\al$ is suitably small, we can construct a suitable weak solution whose error kinetic energy in the strip $\O$ gains one additional order in $\ve$ compared with the corresponding estimate in Section \ref{SEC4}.

We first restate Theorem \ref{thmexistence2} by using Definition \ref{Def2}:

\begin{theorem}[Restatement of Theorem \ref{thmexistence2}]\label{thmexistence2x}
  Assume that \eqref{decayassumption} holds. There exists a universal constant $C_*>0$, independent of $\al$, $\Phi$, $\kappa$, such that when $\e\leq \min\{
    \al^3, \f{\kappa}{12}\}$, and
\bes
\al+\Phi\leq {C}_\ast,
\ees
the system \eqref{nslogerror} with the compatibility condition  \eqref{compatibility1} has a weak solution satisfying
\ba\label{weightedweak}
\|\na (u,v)\cosh(s)\|_{L^2(\O^\e)}^2+\al^{-2}\|(u,v)\cosh s\|_{L^2(\O^\e)}^2\leq C\al^{-4}\kappa^{-1}\Phi^2\e^2.
\ea
Moreover, we derive the following refined trace and kinetic estimates:
\begin{align}
&\|(u,v)(s,\pm\al)\cosh s\|^2_{L^2(\Sigma_\pm)}\leq C \al^{-4} \kappa^{-1}\Phi^2 \e^3, \label{wrefine0x}\\
&\|(u,v)\|^2_{L^2(\O)}\leq C \al^{-6}\Phi^2\left(1+\kappa^{-2}\right)\e^3. \label{wrefine0y}
\end{align}

\end{theorem}

\qed

In the remainder of this section, we complete the proof of Theorem \ref{thmexistence2x}. Using the Banach fixed-point theorem, the existence proof of Theorem \ref{thmexistence2x} follows essentially the same argument as that of Proposition \ref{jh1prop}. Hence, we omit the details and only provide the crucial a priori estimates \eqref{weightedweak} and \eqref{wrefine0x}--\eqref{wrefine0y} in Subsections \ref{Sec5.1} and \ref{Sec5.2}, respectively.

\subsection{The weighted estimate}\l{Sec5.1}

Taking the inner product of \eqref{nslogerror}$_{1,2}$ with $\cosh(2s)(u,v)$, integrating the resulting equations over $\O^\e$, and using integration by parts as in \eqref{042201}, we obtain
\begin{align}
&\int_{\O^\e} \lt( (\p_s u)^2+(\p_s v)^2+ (\p_\th u-v)^2+(\p_\th v+u)^2-2(u^2+v^2) \rt)\cosh(2s)d\th ds\nn\\
&+\int_{\O^\e} e^s\big((u\p_s+v\p_\th) u^0 u+u^0 v^2\big)\cosh(2s)d\th ds-\int_{\O^\e} (u^2+v^2) (u+u^0) e^s\sinh (2s)d\th ds\nn\\
&+\int^{+\i}_{-\i}[\p_\th u]u\Big|_{\Sigma_\pm} \cosh(2s) ds+\int^{+\i}_{-\i}[\p_\th v-e^sp]v\Big|_{\Sigma_\pm} \cosh(2s)ds\nn\\
&-2\int_{\O^\e} p e^s u\sinh(2s)d\th ds=0. \label{wweaks1}
\end{align}

Let $\e\leq\al^3 $. For fixed $\al\in (0,\f{\pi}{100})$ (which is ensured by choosing $\t{C}_*$  sufficiently small), using Minkowski's inequality  and Wirtinger's inequality in Lemma \ref{lemwir}, we have
\bn
\|\p_\th u\|_{L^2_\th}+ \|\p_\th v\|_{L^2_\th}\leq& \|\p_\th u-v\|_{L^2_\th}+\|\p_\th v+u\|_{L^2_\th}+\|u\|_{L^2}+\|v\|_{L^2_\th}\nn\\
                                       \leq &  \|\p_\th u-v\|_{L^2_\th}+\|\p_\th v+u\|_{L^2_\th}+\f{4\al}{\pi}\lt(\|\p_\th u\|_{L^2_\th }+\|\p_\th v\|_{L^2_\th}\rt).
\en
This indicates that
\begin{align}
&\|\p_\th u\|_{L^2_\th}+ \|\p_\th v\|_{L^2_\th}\leq \f{\pi}{\pi-4\al}\lt(\|\p_\th u-v\|_{L^2_\th}+\|\p_\th v+u\|_{L^2_\th}\rt)\nn\\
\leq& 2\lt(\|\p_\th u-v\|_{L^2_\th}+\|\p_\th v+u\|_{L^2_\th}\rt). \label{wweaks2}
\end{align}
In view of Wirtinger's inequality and \eqref{wweaks2}, we have
\begin{align}
\|u\|_{L^2_\th}+ \|v\|_{L^2_\th}\leq \f{8\al}{\pi}\lt(\|\p_\th u-v\|_{L^2_\th}+\|\p_\th v+u\|_{L^2_\th}\rt)\,. \label{wweaks3}
\end{align}
Inserting \eqref{wweaks2} and \eqref{wweaks3} into \eqref{wweaks1}, we can obtain that
\ba\label{wweaks4}
&\|\na (u,v)\cosh s\|^2_{L^2}+ \al^{-2}\|( u, v)\cosh s\|^2_{L^2}\\
\ls &\underbrace{\lt|\int_{\O^\e}\Big(e^s(u\p_s+v\p_\th) u^0 u\cosh(2s) +e^s u^0 v^2\cosh(2s)\Big)d\th ds\rt|}_{J_1}\\
&+\underbrace{\lt|\int_{\O^\e} (u^2+v^2) (u+u^0) e^s\sinh (2s)d\th ds\rt|}_{J_2}\\
&+\underbrace{\lt|\int^{+\i}_{-\i}[\p_\th u]u\Big|_{\Sigma_\pm} \cosh(2s) ds\rt|}_{J_3}+\underbrace{\lt|\int_{\O^\e} p e^s u\sinh(2s)d\th ds\rt|}_{J_4}\,.
\ea

Next, we estimate the terms $J_1$ through $J_4$ one by one. Recalling
\[
u^0=\left\{
\begin{aligned}
&e^{-s}f(\th)\,,&\q\text{for}\q\th\in(-\al,\al)\,;\\
&0\,,&\text{else}\,,
\end{aligned}
\right.
\]
using the estimate in \eqref{jh2ex} and Cauchy's inequality, we have
\begin{align}
|J   _1|\ls \f{\Phi}{\al^{2}}\|(u,v)\cosh s\|^2_{L^2(\O^\e)}\,. \label{wweak5}
\end{align}
By Young's inequality, \eqref{jh2ex}, and the Gagliardo--Nirenberg inequality
\bes
\|f\|^3_{L^3(\O^\e)}\ls \|f\|^3_{H^1(\O^\e)}
\ees
we have

\ba\l{wweak6}
|J_2|\ls& \int_{\O^\e }(|u|^3+|v|^3)\cosh(2s)d\th ds+ \f{\Phi}{\al}\int_{\O^\e }(u^2+v^2)\cosh(2s)d\th ds\\
 \ls&\|\na( u,v)\cosh s\|^3_{L^2}+\|( u, v)\cosh s\|^3_{L^2}+\f{\Phi}{\al}\|( u, v)\cosh s\|^2_{L^2}\,.
\ea
Using the estimate in \eqref{jh2ex}, together with \eqref{compatibility1}, the zero boundary condition, H\"{o}lder's inequality and Cauchy's inequality, we obtain
\bn
|J_3| \ls& {\f{\Phi}{\al^{2}}}\int^{+\i}_{-\i} e^{-s}\Big|\int^{\pm\al\pm\th^\e_{\pm}(s)}_{\pm\al}\p_\th u(s,\th)d\th\Big|\cosh (2s) ds\nn\\
     \ls& {\f{\Phi}{\al^{2}}}\lt|\int^{+\i}_{-\i} e^{3|s|}\lt(\int^{\pm\al\pm\th^\e_{\pm}(s)}_{\pm\al}(\p_\th u(s,\th))^2d\th\rt)^{1/2}|\th^\e_{\pm}(s)|^{1/2} ds\rt|\nn\\
     \ls& {\f{\Phi}{\al^{2}}}\int^{+\i}_{-\i} e^{3|s|}\lt(\int^{\pm\al\pm\th^\e_{\pm}(s)}_{\pm\al}(\p_\th u(s,\th))^2d\th\rt)^{1/2}\e^{1/2} e^{-\f{|s|\kappa}{2\e}} ds\,.\nn\\
\en
Choosing $\e\leq \f{\kappa}{12}$ and applying Cauchy's inequality, we obtain
\ba\label{wweak6ex}
     |J_3|\ls&{\f{\Phi}{\al^{2}}}\sqrt{\e}\int^{+\i}_{-\i} \lt(\int^{\th^\e_{\pm}(s)}_{\pm\al}(\p_\th u(s,\th))^2d\th\rt)^{1/2} e^{-\f{|s|\kappa}{4\e}} ds\\
     \leq& \dl\|\p_\th u\|^2_{L^2(\O^\e)}+C_\dl\e\Phi^2\al^{-4}\int^{+\i}_{-\i}e^{-\f{|s|\kappa}{2\e}} ds\\
     \leq& \dl \|\p_\th u\|^2_{L^2(\O^\e)}+C_\dl {\e^2}\Phi^2\al^{-4}\kappa^{-1}.
\ea

Combining the estimates in \eqref{wweak5}, \eqref{wweak6} and \eqref{wweak6ex}, we obtain that
\begin{align}
  &|J_1|+|J_2|+|J_3|\nn\\
  \ls& \dl \|\p_\th u\cosh(s)\|^2_{L^2(\O^\e)}+C_\dl {\e^2}\Phi^2\al^{-4}\kappa^{-1}\nn\\
         &+\f{\al\Phi+\Phi}{\al^2}\|(u,v)\cosh(s)\|^2_{L^2(\O^\e)}+\|\na(u,v)\cosh(s)\|^3_{L^2(\O^\e)}. \label{wweak6ex1}
\end{align}

In order to estimate $J_4$, we need the following Bogovskii lemma.
{
\begin{lemma}\label{div}
Given $j\in\bZ$, let $\O^\e_j=\{(s,\theta)\in\O^\e\,:\,j\al<s<(j+1)\al\}$ and $h\in L^2(\O^\e_j)$ with
\bn
\int_{\O^\e_j}h d s d\theta=0\,.
\en
Then there exists a vector-valued function $\bl{\varphi}_j=(\phi_j,\psi_j)\,:\,\O^\e_j\to\mathbb{R}^2$ belonging to $H^1_0(\O^\e_j)$ such that
\bn
&\p_s\phi_j+\p_\theta\psi_j=h\,,
\en
and
\ba\l{vf2}
\|\na\bl{\varphi}_j\|_{L^2(\O^\e_j)}\leq \frak{C}_\ast\|h\|_{L^2(\O^\e_j)}\,,\q\|\bl{\varphi}_j\|_{L^2(\O^\e_j)}\leq  \al\frak{C}_\ast \|h\|_{L^2(\O^\e_j)}\,,
\ea
where the Bogovskii constant $\frak{C}_\ast$ is independent of $\e$, $\al$ and $j$.
\end{lemma}

\begin{proof}
The proof is a direct application of Theorem III.3.1 of \cite[Section III.3]{Galdi2011} and Lemma \ref{lemwir}. We make a brief explanation here. From Theorem III.3.1 of \cite[Section III.3]{Galdi2011}, there exists $\bl{\varphi}_j$ satisfying
\bn
\p_s\phi_j+\p_\theta\psi_j=h,\text{ and }\q  \|\na\bl{\varphi}\|_{L^2(\O^\e_j)}\leq C_\ast\|h\|_{L^2(\O^\e_j)},
\en
where
\be\label{bogproof1}
C_\ast=C\lt(\f{\dl({\O^\e_j})}{R}\rt)^2\lt(1+\f{\dl({\O^\e_j})}{R}\rt).
\ee
The constant $C$ depends on the Lipschitz norm of the boundary graph, $R$ denotes the radius of a ball contained in $\O^\e_j$, and $\dl({\O^\e_j})$ is the diameter of the domain. For our domain ${\O^\e_j}$, since
\man{
|\th^\e_\pm(s)|\ls \e\leq 1, \label{bogproof2}
}
and by applying \eqref{decayassumption}, we have
\man{
|(\th^\e_\pm)^\prime(s)|=e^{\f{|s|}{\e}}|\t{\gamma}'_{\pm}(e^{\f{|s|}{\e}})|\leq M_1. \label{bogproof3}
}
From this we see that $C$ in \eqref{bogproof1} is uniformly bounded and independent of $j$ and $\e$. Besides, the maximum radius of a ball contained in ${\O^\e_j}$ we can choose is approximately $\al$, so
\be
R\approx \al. \label{bogproof4}
\ee
Also, it is clear that
\be
\dl({\O^\e_j})\approx \al.  \label{bogproof5}
\ee
Inserting \eqref{bogproof2}, \eqref{bogproof3}, \eqref{bogproof4} and \eqref{bogproof5} into \eqref{bogproof1}, we see that
\bn
C_\ast\leq C_{{M_1}} (1+\f{\al}{\al})\ed C_{\text{uni}}.
\en
Here the constant $C_{\text{uni}}$ is independent of $\e$, $\al$ and $j$. Finally, the $L^2$ norm estimate of $\bl{\varphi}$ follows from Lemma \ref{lemwir}, which indicates
\bes
 \|\bl{\varphi}\|_{L^2(\O^\e_j)} \leq C_{\text{Wri}}\al\|\p_s\bl{\varphi}\|_{L^2(\O^\e_j)}\leq C_{\text{Wri}}C_{\text{uni}}\al\|h\|_{L^2(\O^\e_j)}.
\ees
Hence by denoting it by $\frak{C}_\ast\ed C_{\text{uni}}(1+C_{\text{Wri}})$, we obtain \eqref{vf2}.
\end{proof}
}

Since $u$ has zero flux, the integral of $(e^{3s}-e^{-s})u$ over each $\theta-$section vanishes, and hence its integral over the entire domain $\O^\e_j$ is zero $\forall j\in\bZ$. Using Lemma \ref{div}, there exists a sequence of vector fields $\{\bl{\varphi}_j\}_{j\in\mathbb{Z}}$ such that:
\be\label{Auxp}
\left\{
\begin{aligned}
&\p_s\phi_{j}+\p_\theta\psi_{j}=(e^{3s}-e^{-s})u\,,\q\text{in }\O^\e_j\,,\\
&\bl{\varphi}_j\in H_0^1(\O^\e_j)\,,\\
&\|\bl{\varphi}_j\|_{H^1(\O^\e_j)}\leq\frak{C}_\ast(e^{3j\al}+e^{-j\al})\|u\|_{L^2(\O^\e_j)}\,.
\end{aligned}
\right.
\ee
 Thus we decompose $\O^\e=\cup_{j\in\bZ}\O^\e_j$ and estimate $J_4$ as follows. Using \eqref{Auxp}, integration by parts and then replacing $\na p$ with the system \eqref{nslogerror}, we obtain
\[
\begin{split}
|J_4|\ls&\Big|\sum_{j\in\bZ}\int_{\O^\e_j} pu(e^{3s}-e^{-s})dsd\th\Big|\\
=&\Big|\sum_{j\in\bZ}\int_{\O^\e_j}(\p_s\phi_{j}+\p_\theta \psi_{j})p dsd\th\Big|=\Big|\sum_{j\in\bZ}\int_{\O^\e_j}\phi_{j}\p_sp+\psi_{j}\p_\theta p dsd\th\Big|\\
\leq &\Big|\sum_{j\in\bZ}\int_{\O^\e_j}\phi_{j}\left\{e^{-s}\left(\Dl_s u-2\p_\th v-u\right)-\lt[({u}+u^0)\p_s u+{v}\p_\th u+ (u\p_s+v\p_\th) u^0-v^2\rt]\rt\}\Big|\\
&+\Big|\sum_{j\in\bZ}\int_{\O^\e_j}\psi_{j}\left\{e^{-s}\left(\Dl_s v+2\p_\th u-v\right)-\lt[({u}+u^0)\p_s v+{v}\p_\th v+(u+u^0)v\rt]\rt\}\Big|.
\end{split}
\]

Using integration by parts for terms $\phi_j e^{-s} \Dl_s u$, $\psi_j e^{-s} \Dl_s v$ and $u^0\p_s u\phi_j$, $u^0\p_s v\psi_j$, and then noticing that
\[
u^0=\left\{
\begin{aligned}
&e^{-s}f(\th)\,,&\q\text{for}\q\th\in(-\al,\al)\,;\\
&0\,,&\text{else}\,,
\end{aligned}
\right.
 \]
we can obtain the following estimates by applying Proposition \ref{jh1prop}:
\begin{align}
|J_4|\ls& \underbrace{\sum_{j\in\bZ} \int_{\O^\e_j}\lt(|\na u|+|\na v|+|u|+|v|\rt)(|\na\bl{\varphi}_j|+|\bl{\varphi}_j|)e^{-s}d\th ds}_{J_{41}}\nn\\
&+\underbrace{\sum_{j\in\bZ} \int_{\O^\e_j}\lt(u^2+v^2+|\na (u,v)||(u,v)|\rt)|\bl{\varphi}_j|d\th ds}_{J_{42}} \label{wweak7} \\
&+\underbrace{\f{\Phi}{\al^2}\sum_{j\in\bZ} \int_{\O^\e_j}|(u,v)|\cdot\lt|(\na \bl{\varphi}_j,\bl{\varphi}_j)\rt|e^{-s}d\th ds.}_{J_{43}}\nn
\end{align}
Since in $\O^\e_j$, $e^{\pm s}\simeq e^{\pm j\al}$.  By H\"{o}lder's inequality and \eqref{Auxp}, we have
\begin{align}
|J_{41}|\ls& \sum_{j\in\bZ}\lt(\|\na(u,v)\cosh s\|_{L^2(\O^\e_j)}+\|(u,v)\cosh s\|_{L^2(\O^\e_j)}\rt)\nn\\
&\cdot\lt(\|\na\bl{\varphi}_j\f{e^{-s}}{\cosh s}\|_{L^2(\O^\e_j)}+\|\bl{\varphi}_j\f{e^{-s}}{\cosh s}\|_{L^2(\O^\e_j)}\rt)\nn\\
\ls&  \sum_{j\in\bZ}\lt(\|\na(u,v)\cosh s\|_{L^2(\O^\e_j)}+\|(u,v)\cosh s\|_{L^2(\O^\e_j)}\rt)\|\bl{\varphi}_j\|_{H^1(\O^\e_j)}(e^{2j\al}+1)^{-1} \nn\\
\ls& {\frak{C}_\ast}\sum_{j\in\bZ}\lt(\|\na(u,v)\cosh s\|_{L^2(\O^\e_j)}+\|(u,v)\cosh s\|_{L^2(\O^\e_j)}\rt)\|u\|_{L^2(\O^\e_j)}{(e^{j\al}+e^{-j\al})}\nn\\
\ls& {\frak{C}_\ast} \sum_{j\in\bZ}\lt(\|\na(u,v)\cosh s\|_{L^2(\O^\e_j)}+\|(u,v)\cosh s\|_{L^2(\O^\e_j)}\rt)\|u\cosh s\|_{L^2(\O^\e_j)}.\nn
\end{align}
Then, by Cauchy's inequality, we have the following estimate.
\begin{align}
|J_{41}|\ls&  \dl \sum_{j\in\bZ}\|\na(u,v)\cosh s\|^2_{L^2(\O^\e_j)}+C_\dl {\frak{C}^2_\ast}\sum_{j\in\bZ}\|(u,v)\cosh s\|^2_{L^2(\O^\e_j)}\nn \\
\ls&\dl \|\na(u,v)\cosh s\|^2_{L^2(\O^\e)}+C_\dl {\frak{C}^2_\ast} \|(u,v)\cosh s\|^2_{L^2(\O^\e)}. \l{wweak8}
\end{align}
Similarly to \eqref{wweak8}, we infer from H\"{o}lder's inequality and \eqref{Auxp} that
\begin{align}
|J_{43}|\ls& \f{\Phi}{\al^2}\sum_{j\in\bZ}\|(u,v)\cosh s\|_{L^2(\O^\e_j)}\lt(\|\na\bl{\varphi}_j\f{e^{-s}}{\cosh s}\|_{L^2(\O^\e_j)}+\|\bl{\varphi}_j\f{e^{-s}}{\cosh s}\|_{L^2(\O^\e_j)}\rt)\nn\\
\ls& \f{\Phi}{\al^2}\sum_{j\in\bZ}\|(u,v)\cosh s\|_{L^2(\O^\e_j)}\|\bl{\varphi}_j\|_{H^1(\O^\e_j)}(e^{2j\al}+1)^{-1} \l{wweak9} \\
\ls& {\frak{C}_\ast}\f{\Phi}{\al^2}\sum_{j\in\bZ}\|(u,v)\cosh s\|_{L^2(\O^\e_j)}\|u\|_{L^2(\O^\e_j)}(e^{j\al}+e^{-j\al})\nn\\
\ls&   \f{\Phi}{\al^2}\frak{C}_\ast\sum_{j\in\bZ}\|(u,v)\cosh s\|^2_{L^2(\O^\e_j)}\ls{\f{\Phi}{\al^2}\frak{C}_\ast}\|(u,v)\cosh s\|^2_{L^2(\O^\e)}.\nn
\end{align}
Finally, by H\"{o}lder's inequality, Sobolev embedding, \eqref{Auxp} and Young's inequality, we have
\begin{align}
|J_{42}|\ls& \sum_{j\in\bZ}\|(u,v)\|^2_{L^4(\O^\e_j)}\|\bl{\varphi}_j\|_{L^2(\O^\e_j)}+\|(u,v)\|_{L^4(\O^\e_j)}\|\na(u,v)\|_{L^2(\O^\e_j)}\|\bl{\varphi}_j\|_{L^4(\O^\e_j)}\nn\\
       \ls& \sum_{j\in\bZ}\lt(\|(u,v)\|_{L^2(\O^\e_j)}\|\na(u,v)\|_{L^2(\O^\e_j)}\|\bl{\varphi}_j\|_{L^2(\O^\e_j)}\rt)\nn\\
       &+ \|(u,v)\|^{1/2}_{L^2(\O^\e_j)}\|\na(u,v)\|^{3/2}_{L^2(\O^\e_j)}\|\bl{\varphi}_j\|^{1/2}_{L^2(\O^\e_j)}\|\na\bl{\varphi}_j\|^{1/2}_{L^2(\O^\e_j)}\nn\\
       \ls& \frak{C}_\ast \sum_{j\in\bZ}\lt(\|(u,v)\|^2_{L^2(\O^\e_j)}\|\na(u,v)\|_{L^2(\O^\e_j)}+ \|(u,v)\|^{3/2}_{L^2(\O^\e_j)}\|\na(u,v)\|^{3/2}_{L^2(\O^\e_j)}\rt)(e^{3j\al}+e^{-j\al})\nn\\
       \ls &  \frak{C}_\ast \|\na(u,v)\cosh s\|^3_{L^2(\O^\e)}+ \frak{C}_\ast\|(u,v)\cosh s\|^3_{L^2(\O^\e)}.\label{wweak10}
\end{align}
Inserting \eqref{wweak8},\eqref{wweak9} and \eqref{wweak10} into \eqref{wweak7}, we see that
\begin{align}
|J_4|\ls& \dl \|\na(u,v)\cosh s\|^2_{L^2(\O^\e)}+ \frak{C}_\ast\|\na(u,v)\cosh s\|^3_{L^2(\O^\e)}+ \frak{C}_\ast\|(u,v)\cosh s\|^3_{L^2(\O^\e)}\nn\\
      &+\lt(C_\dl{\frak{C}^2_\ast}+\f{\Phi}{\al^2}\frak{C}_\ast\rt)\|(u,v)\cosh s\|^2_{L^2(\O^\e)}. \label{wweak11}
\end{align}
Combining the estimates in \eqref{wweak6ex1} and \eqref{wweak11}, we have the following estimate.
\begin{align}
&|J_1|+|J_2|+|J_3|+|J_4|\nn\\
\leq& C\dl\|\na(u,v)\cosh s\|^2_{L^2(\O^\e)}\nn\\
    &+C\lt(\frak{C}_\ast+1\rt)\|\na(u,v)\cosh s\|^3_{L^2(\O^\e)}+C{\frak{C}_\ast}\|(u,v)\cosh s\|^3_{L^2(\O^\e)} \l{wweak11a}\\
&+\lt(\lt(C_\dl\frak{C}^2_\ast+\f{\Phi}{\al^2}\frak{C}_\ast\rt)+\f{\al\Phi+\Phi}{\al^2}\rt)\|(u,v)\cosh s\|^2_{L^2(\O^\e)}+C_\dl {\e^2}\Phi^2\al^{-4}\kappa^{-1}\,. \nn
\end{align}
Inserting \eqref{wweak11a} into \eqref{wweaks4} and choosing $\dl$ sufficiently small, we obtain the existence of a sufficiently small constant ${C}_\ast$ such that, whenever
\be
\al+\Phi\leq {C}_\ast\,,  \label{wweak11b}
\ee
we have
\begin{align}
&\|\na( u, v)\cosh s\|^2_{L^2}+\al^{-2}\|( u, v)\cosh s\|^2_{L^2}\nn\\
\ls& \|\na( u, v)\cosh s\|^3_{L^2}+\|( u, v)\cosh s\|^3_{L^2}+\al^{-4}{\e^2}\Phi^2\kappa^{-1}. \nn
\end{align}
This a priori estimate indicates that under the assumption \eqref{wweak11b}, the system \eqref{nslogerror} has a weak solution satisfying
\be
\|\na( u, v)\cosh s\|^2_{L^2}+\al^{-2}\|( u, v)\cosh s\|^2_{L^2}\ls \al^{-4}\e^2\Phi^2\kappa^{-1}. \label{wweak12}
\ee
Indeed, we construct the solution by a standard Banach fixed point argument in the weighted space
with norm
\[
\|(u,v)\|_X^2
:=
\|\nabla(u,v)\cosh s\|_{L^2(\Omega^\varepsilon)}^2
+\alpha^{-2}\|(u,v)\cosh s\|_{L^2(\Omega^\varepsilon)}^2 .
\]
For $\alpha+\Phi$ sufficiently small, the corresponding nonlinear map is a contraction
on the ball
\[
\|(u,v)\|_X\le C\alpha^{-2}\Phi\varepsilon\kappa^{-1/2}.
\]
Hence the fixed point obtained in this ball satisfies
\[
\|\nabla(u,v)\cosh s\|_{L^2(\Omega^\varepsilon)}^2
+\alpha^{-2}\|(u,v)\cosh s\|_{L^2(\Omega^\varepsilon)}^2
\le
C\alpha^{-4}\Phi^2\varepsilon^2\kappa^{-1}.
\]
We omit the details of this standard contraction argument. This completes the proof of \eqref{weightedweak}.

\subsection{The refined estimate}\l{Sec5.2}
In the following, we prove the estimates  \eqref{wrefine0x} and \eqref{wrefine0y} based on the weighted estimate \eqref{wweak12}.  The fundamental theorem of calculus is applied to convert the line integral into a volume integral, then the trace estimate \eqref{wrefine0x} follows directly from the boundary condition. The key idea in the proof of \eqref{wrefine0y} is to use the refined boundary trace estimate \eqref{wrefine0x}, which yields an extra power of $\varepsilon$ in the boundary error. This improved boundary control is then converted, through an adjoint argument, into a stronger $L^2$ estimate for the error in the bulk domain.

\subsubsection{The improved boundary estimate}
First, we derive a weighted trace estimate on the line $\Sigma_\pm$. Using Poincar\'{e} inequality in $\theta$ direction, it follows that
\begin{align}
\big|\int^{\pm \al\pm\th^\e_\pm(s)}_{\pm \al}(u^2+v^2)d\th\big|\les \e^2\big|\int^{\pm \al\pm\th^\e_\pm(s)}_{\pm \al}\big((\p_\th u)^2+(\p_\th v)^2\big)d\th\big|. \label{poincare}
\end{align}
Since $(u,v)=\bm{0}$ on the rough boundaries, by the boundary condition, H\"{o}lder inequality and \eqref{poincare}, we have
\ba\label{wrefine}
\|(u,v)(s,\pm\al)\cosh s\|^2_{L^2(\Sigma_\pm)}=&\lt|\int_{\bR} \int^{\pm\al\pm\th^\e_\pm(s)}_{\th^\e_\pm(s)}\p_\th (u^2+v^2) d\th\cosh (2s) ds\rt|\\
\ls & \lt(\int_{\bR}\int^{\pm\al\pm\th^\e_\pm(s)}_{\th^\e_\pm(s)}(u^2+v^2)\cosh (2s)d\th ds\rt)^{1/2}\\
&\cdot\lt(\int_{\bR}\int^{\pm\al\pm\th^\e_\pm(s)}_{\th^\e_\pm(s)}((\p_\th u)^2+(\p_\th v)^2)\cosh (2s)d\th ds\rt)^{1/2}\\
\ls&\, \e \|\p_\th( u, v)\cosh s\|^{2}_{L^2}.
\ea
Substituting \eqref{wweak12} in \eqref{wrefine}, we obtain that
\begin{align}
\|(u,v)(s,\pm\al)\cosh s\|^2_{L^2(\Sigma_\pm)}\ls \al^{-4}\kappa^{-1}\Phi^2\e^3\,. \label{wrefine0}
\end{align}
This proves \eqref{wrefine0x}.
\subsubsection{The adjoint approach}
Motivated by \cite{BG2008CPAM}, we use the boundary estimate \eqref{wrefine0} and the adjoint method to prove \eqref{wrefine0y}. Let $(w,z,r)$ denote a solution of the following boundary value problem

\be\label{nslogerror2}
\left\{
\begin{aligned}
&e^s\p_s r- \left(\Dl_s w-2\p_\th z-w\right)=u, \q \text{in}\q \O,\\
&e^s\p_\th r- \left(\Dl_s z+2 \p_\th w-{z}\right)=v,\q \text{in}\q \O,\\
&\p_s (e^sw)+\p_\th (e^sz)=0,  \q \text{in}\q \O,\\
&\lt(w,z\rt)=0, \q \text{on}\q \p\O.
\end{aligned}
\right.
\ee
We now derive an energy estimate for the system \eqref{nslogerror2}. Multiplying \eqref{nslogerror2}$_{1}$ and \eqref{nslogerror2}$_{2}$ by $w$ and $z$, respectively, integrating the resulting equations over $\O$, and using integration by parts, we obtain
\bn
&\int_{\O} \lt( (\p_s w)^2+(\p_s z)^2+ (\p_\th w-z)^2+(\p_\th z+w)^2 \rt)d\th ds=\int_\O (uw+vz)d\th ds\,.
\en
Using Lemma \ref{LemCP} in the same way as in the derivation of \eqref{wweaks4}, we obtain
\bn
\|\na( w, z)\|^2_{L^2(\O)}+\al^{-2}\|( w, z)\|^2_{L^2(\O)}\ls \big|\int_\O (uw+vz)d\th ds\big|.
\en
Then, by Cauchy's inequality, we obtain
\be
\|\na( w, z)\|^2_{L^2(\O)}+\al^{-2}\|( w, z)\|^2_{L^2(\O)}\ls \al^2\|( u, v)\|^2_{L^2(\O)}. \label{wrefine1}
\ee
%
%
%

Next we derive the second-order derivative estimates for $(w,z)$. For convenience, we rewrite \eqref{nslogerror2} as follows:
\be\label{wrefine5}
\left\{
\begin{array}{l}
\p_s r- \Dl_s(e^{-s}w)=2\p_s(e^{-s}w)-2e^{-s}\p_\th z+e^{-s}u, \q \text{in}\q \O,\\
\p_\th r- \Dl_s(e^{-s}z)=2\p_s(e^{-s}z)+2e^{-s}\p_\th w+e^{-s}v,\hskip.2cm\q \text{in}\q \O,\\
\p_s (e^{-s}w)+\p_\th (e^{-s}z)=-2e^{-s}w,  \q \text{in}\q \O,\\
\lt(w,z\rt)=0, \q \text{on}\q \p\O.
\end{array}
\right.
\ee

Our strategy for obtaining the second-order derivative estimates for \eqref{wrefine5} is first to cut off the domain into $\cup_{j\in\bZ} \O_j$ with { $\O_j:=[j\al,(j+1)\al]\times[-\al,\al]$}, then to apply estimates for the Stokes system on uniformly locally bounded domains. A simple gluing argument then gives the $H^2$ estimate for $(w,z)$.

To localize the problem \eqref{wrefine5}, we introduce the following sequence of cut-off functions $\chi_j(s)\in C^\i_c(\bR)$, which satisfies
{
\bn
\chi_j(s)=\lt\{
\bali
&1,\q s\in [j\al,(j+1)\al]\,;\\
&0,\q s\leq (j-1)\al \q\text{ or }\q s>(j+2)\al\,,
\eali
\rt.
\en
and
\be\label{cutoff}
|\chi_j^\prime|\ls \f{1}{\al},\q |\chi^{\prime\prime}_j|\ls \f{1}{\al^2}.
\ee
}
We also let $c_j$ denote the mean value of $r$ in the cut-off domain $\t{\O}_j:=[(j-1)\al,(j+2)\al]\times[-\al,+\al]$. From \eqref{wrefine5}, we have
\be\l{STOKES0428}
\lt\{
\begin{aligned}
&\p_s (\chi_j(r-c_j))- \Dl_s(e^{-s}w\chi_j)=\frak{f}_j\,,\\
&\p_\th (\chi_j(r-c_j))- \Dl_s(e^{-s}z\chi_j)=\frak{g}_j\,,\\
&\p_s (e^{-s}w\chi_j)+\p_\th (e^{-s}z\chi_j)=\frak{h}_j\,, \\
&\lt(e^{-s}\chi_jw,e^{-s}\chi_jz\rt)\big|_{\p\t{\O}_j}=0\,.
\end{aligned}
\right.
\ee
Here, by direct calculation, $\frak{f}_j$,  $\frak{g}_j$ and $\frak{h}_j$ are given as follows:
\begin{align}
&\frak{f}_j:=\chi_j\lt(2(e^{-s}w)_s-2e^{-s}z_\th\rt)+\chi^\prime_j(r-c_j)-2\chi^\prime_j(e^{-s}w)_s-\chi^{\prime\prime}_je^{-s}w+\chi_j e^{-s}u\,,\nn\\
&\frak{g}_j:=\chi_j\lt(2(e^{-s}z)_s+2e^{-s}w_\th\rt)-2\chi^\prime_j(e^{-s}z)_s-\chi^{\prime\prime}_je^{-s}z+\chi_j e^{-s}v\,,\l{RHS0428}\\
&\frak{h}_j:=-2e^{-s}w\chi_j+\chi^\prime_je^{-s}w\,.\nn
\end{align}
For the problem \eqref{STOKES0428}, we use the standard $H^2$ estimate in $\t{\O}_j$  for stationary Stokes equations to get
{ \ba
\|e^{-s}\chi_j(w,z)\|^2_{\dot{H}^2(\t{\O}_j)}\leq C\left(\|(\frak{f}_j,\frak{g}_j)\|^2_{L^2(\t{\O}_j)}+\al^{-2}\|\frak{h}_j\|^2_{L^2(\t{\O}_j)}+\|\na \frak{h}_j\|^2_{L^2(\t{\O}_j)}\right).\label{wrefine6}
\ea
The factor $\al^{-2}$ can be obtained by a scaling argument by scaling the domain $\t{\O}_j$ to the standard rectangle with width $3$ and height $2$.
}
Here note that the constant $C$ is independent of $j$ since all $\t{\O}_j$ are congruent. A detailed proof of \eqref{wrefine6} can be found in \cite[Chapter IV]{Galdi2011}. Meanwhile, in view of \eqref{RHS0428} and \eqref{cutoff}, we derive by direct calculations that
{
\begin{align}\label{wrefine6x}
\|(\frak{f}_j,\frak{g}_j)\|^2_{L^2(\t{\O}_j)}\ls& e^{-2j\al}\lt(\al^{-2}\|(\na w,\na z)\|^2_{L^2(\t{\O}_j)}+\al^{-4}\|(w,z)\|^2_{L^2(\t{\O}_j)}+\|(u,v)\|^2_{L^2(\t{\O}_j)}\rt)\nn\\
                                                &+\al^{-2}\|r-c_j\|^2_{L^2(\t{\O}_j)},\nn\\
\|\frak{h}_j\|^2_{L^2(\t{\O}_j)}\ls& e^{-2j\al}\al^{-2}\|w\|^2_{L^2(\t{\O}_j)},\nn\\
\|\na \frak{h}_j\|^2_{L^2(\t{\O}_j)}\ls& e^{-2j\al}\al^{-2}\|\na w\|^2_{L^2(\t{\O}_j)}+e^{-2j\al}\al^{-4}\| w\|^2_{L^2(\t{\O}_j)}.
\end{align}
}
Since $e^{\pm s}\simeq e^{\pm j\al}$ in $\t{\O}_j$, inserting \eqref{wrefine6x} into the right-hand side of \eqref{wrefine6}, we show that
{
\begin{align}
\|(w,z)\|^2_{\dot{H}^2({\O}_j)}\ls&\al^{-2}\|(\na w,\na z)\|^2_{L^2(\t{\O}_j)}+\al^{-4}\|(w,z)\|^2_{L^2(\t{\O}_j)}+\|(u,v)\|^2_{L^2(\t{\O}_j)}\nn\\
                         &+e^{2j\al}\al^{-2}\|r-c_j\|^2_{L^2(\t{\O}_j)}.\label{wrefine7}
\end{align}
}
Now it remains to estimate $r-c_j$ in $\t{\O}_j$. From Lemma \ref{div}, there exists $\bl{\varphi}_j=(\phi_j,\psi_j)\in H_0^1(\tilde{\O}_j)$ such that
\be\l{EE0428}
\p_s \phi_j+\p_\th\psi_j=r-c_j,\q \|\na\bl{\varphi}_j\|_{L^2(\t{\O}_j)}+\al^{-1}\|\bl{\varphi}_j\|_{L^2(\t{\O}_j)}\ls \|r-c_j\|_{L^2(\t{\O}_j)}.
\ee From \eqref{wrefine5}, integration by parts, Lemma \ref{div} and H\"{o}lder's inequality, we have
\bn
\|r-c_j\|^2_{L^2(\t{\O}_j)}=&-\int_{\t{\O}_j} \na(r-c_j)\cdot\bl{\varphi}_jd\th ds\\
=&\int_{\t{\O}_j} \lt[- \Dl_s(e^{-s}w)-2\p_s(e^{-s}w)+2e^{-s}z_\th-e^{-s}u\rt]\phi_j d\th ds\\
&+\int_{\t{\O}_j}\lt[-\Dl_s(e^{-s}z)-2\p_s(e^{-s}z)-2e^{-s}w_\th-e^{-s}v\rt]\psi_jd\th ds\\
=& \int_{\t{\O}_j} \lt[-2\p_s(e^{-s}w)+2e^{-s}z_\th-e^{-s}u\rt]\phi_j+\lt[-2\p_s(e^{-s}z)-2e^{-s}w_\th-e^{-s}v\rt]\psi_jd\th ds\\
 &+\int_{\t{\O}_j}\na(e^{-s}w)\cdot\na\phi_j+\na(e^{-s}z)\cdot\na\psi_jd\th ds\\
\ls& {e^{-j\al}\|\bl{\varphi}_j\|_{H^1(\t{\O}_j)}\|(\na w,\na z,w,z) \|_{L^2(\t{\O}_j)}+e^{-j\al}\|\bl{\varphi}_j\|_{L^2(\t{\O}_j)}\|(u,v) \|_{L^2(\t{\O}_j)}}\,.
\en
This, together with \eqref{EE0428}, indicates that
\be
{\|r-c_j\|_{L^2(\t{\O}_j)}\ls e^{-j\al}\|(\na w,\na z,w,z) \|_{L^2(\t{\O}_j)}+e^{-j\al}\al\|(u,v) \|_{L^2(\t{\O}_j)}}.\label{wrefine8}
\ee
Thus by inserting \eqref{wrefine8} into \eqref{wrefine7}, we obtain
\be
\|(w,z)\|^2_{\dot{H}^2({\O}_j)}\ls {\al^{-2}\|(\na w,\na z)\|^2_{L^2(\t{\O}_j)}+\al^{-4}\|(w,z)\|^2_{L^2(\t{\O}_j)}+\|(u,v)\|^2_{L^2(\t{\O}_j)}}. \label{wrefine9}
\ee
Summing over $j\in \bZ$, estimate \eqref{wrefine9} indicates
\bn
\|\na^2(w,z)\|^2_{L^2({\O})}\ls{\al^{-2}\|(\na w,\na z)\|^2_{L^2(\O)}+\al^{-4}\|(w,z)\|^2_{L^2(\O)}+\|(u,v)\|^2_{L^2(\O)}}.\label{wrefine10}
\en
This, together with \eqref{wrefine1}, indicates
\begin{align}
\|\na^2(w,z)\|^2_{L^2({\O})}\ls \|( u, v)\|^2_{L^2(\O)}. \label{wrefine11}
\end{align}

\subsubsection{Completion of the proof}

Finally, we use the adjoint system \eqref{nslogerror2} and estimates \eqref{wrefine1}, \eqref{wrefine11} to prove the estimate \eqref{wrefine0y}.  Directly from the system \eqref{nslogerror2}, after integration by parts, we have
\begin{align*}
&\int_\O (|u|^2+|v|^2)d\th ds\\
=&\int_\O \Big(u[e^s\p_s r- \left(\Dl_s w-2\p_\th z-w\rt)]+v[e^s\p_\th r- \left(\Dl_s z+2 \p_\th w-{z}\rt)]\Big)d\th ds\\
  =&\int_\O \big(- w\Dl_s u-2z\p_\th u+uw -z\Dl_s v+2 w\p_\th v+ vz\big)d\th ds\\
  &+ \int_{\p\O}\big( w\p_n u-u\p_n w \big)dl+ \int_{\p\O}\big( z\p_n v-v\p_n z \big)dl\\
  =&\int_\O \Big(-(\Dl_s u-2\p_\th v-u) w -(\Dl_s v +2\p_\th u - v)z \Big)d\th ds\\
  &-\int_{\bR} u\p_\th w\Big|_{\th=-\al}^\al ds-\int_{\bR}v \p_\th z\Big|_{\th=-\al}^\al ds.
\end{align*}
 Then by substituting \eqref{nslogerror}$_{1,2}$ into the first term on the right-hand side of the above equality, and applying the boundary condition $(w,z)\big|_{\th=\pm\al}=0$ together with \eqref{nslogerror2}$_3$, the pressure term can be canceled through integration by parts, yielding
\begin{align}
&\int_\O (|u|^2+|v|^2)\nn\\
   =&\underbrace{-\int_\O e^{s}\lt[({u}+u^0)\p_s u+({v}+v^0)\p_\th u+ (u\p_s+v\p_\th) u^0-v^2\rt]w\,d\th ds }_{K_1}\nn\\
   &\underbrace{-\int_\O e^s\lt[({u}+u^0)\p_s v+({v}+v^0)\p_\th v+ (u\p_s+v\p_\th) v^0 +(u+u^0) v\rt]z\,d\th ds}_{K_2}\label{wrefine13}\\
   &\underbrace{-\int_{\bR} u\p_\th w\Big|_{\th=-\al}^\al ds-\int_{\bR}v \p_\th z\Big|_{\th=-\al}^\al ds}_{K_3} \,.\nn
\end{align}
 We now estimate the terms $K_1$, $K_2$ and $K_3$.
Using integration by parts, H\"{o}lder inequality, and the property of $u_0$ in Proposition \ref{jh1prop}, we see that
\begin{align*}
|K_1|=&\big|\int_\O \big(u\lt(e^s(u+u_0)\p_s w+e^s v \p_\th w\rt)+uwf(\th)-vwf'(\th)+e^sv^2w\big)\,d\th ds\big|\\
    \ls&\|e^{\f{s}{2}}(u,v)\|^2_{L^4}\|(\p_sw,\p_\th w,w)\|_{L^2}+\f{\Phi}{\al^2}\|(u,v)\|_{L^2}\|w\|_{L^2}+\f{\Phi}{\al}\|(u,v)\|_{L^2}\|\p_sw\|_{L^2}\,.
\end{align*}
 By using \eqref{wrefine1}, \eqref{wrefine11}, Sobolev embedding and Cauchy inequality, we have
\bn
    |K_1|\ls& \|(u,v)\|_{L^2}\|e^{\f{s}{2}}(u,v)\|^2_{L^4}+\Phi\|(u,v)\|^2_{L^2} \text{}\nn\\
    \ls& \|(u,v)\|_{L^2}\|\cosh s(u,v)\|_{L^2}\|\cosh s\na(u,v)\|_{L^2}+\Phi\|(u,v)\|^2_{L^2}\nn\\
    \ls& (\dl+\Phi)\|(u,v)\|^2_{L^2}+C\dl^{-1} \|\cosh s(u,v)\|^2_{L^2}\|\cosh s\na(u,v)\|^2_{L^2}\,.\\
    \en
Therefore, in view of \eqref{wweak12}, one deduces
{\ba\l{wrefine14}
    |K_1|\ls&(\dl+\Phi)\|(u,v)\|^2_{L^2}+C\dl^{-1} \al^{-6}\Phi^4 \e^4\kappa^{-2}\,.
\ea}
Akin to $K_1$, by H\"{o}lder's inequality, \eqref{wrefine1}, \eqref{wrefine11}, and Sobolev embedding, we have
 \begin{align*}
|K_2|=&\big|\int_\O \big(v\lt(e^s(u+u_0)\p_s z+e^sv \p_\th z \rt)-e^suvz-f(\th)vz\big)\,d\th ds\big|\\
    \ls&\|e^{\f{s}{2}}(u,v)\|^2_{L^4}\|(\p_sz,\p_\th z,z)\|_{L^2}+\f{\Phi}{\al}\|(u,v)\|_{L^2}\|(\p_s z,z)\|_{L^2}\\
    \ls& \|(u,v)\|_{L^2}\|e^{\f{s}{2}}(u,v)\|^2_{L^4}+\Phi\|(u,v)\|^2_{L^2}.
\end{align*}
Then, by Cauchy's inequality and \eqref{wweak12}, we have the following estimate for $K_2$:
\begin{align}
  |K_2|  \ls& \|(u,v)\|_{L^2}\|\cosh s(u,v)\|_{L^2}\|\cosh s\na(u,v)\|_{L^2}+\Phi\|(u,v)\|^2_{L^2}\nn\\
    \ls& (\dl+\Phi)\|(u,v)\|^2_{L^2}+C\dl^{-1} \|\cosh s(u,v)\|^2_{L^2}\|\cosh s\na(u,v)\|^2_{L^2}\l{wrefine15}\\
    \ls&(\dl+\Phi)\|(u,v)\|^2_{L^2}+C\dl^{-1} \al^{-6}\Phi^4{\e^4}\kappa^{-2}\,.\nn
\end{align}
Finally, using the Cauchy-Schwarz inequality and the trace inequality, one has
\bn
|K_3|\leq& \|(u,v)\|_{L^2(\Sigma_\pm)}\|(\p_\th w,\p_\th z)\|_{L^2(\Sigma_\pm)}\ls \|(u,v)\|_{L^2(\Sigma_\pm)}\|(w,z)\|_{H^2(\O)}\,.
\en
Thus one deduces by using Cauchy's inequality, \eqref{wrefine0}, \eqref{wrefine1} and \eqref{wrefine11} that
\ba\l{wrefine16}
|K_3|\leq&\|(u,v)\|_{L^2(\Sigma_\pm)}\|(u,v)\|_{L^2(\O)}\\
\leq & \dl\|(u,v)\|^2_{L^2(\O)}+C\dl^{-1}\|(u,v)\|^2_{L^2(\Sigma_\pm)}\\
\ls& \dl\|(u,v)\|^2_{L^2(\O)}+C\dl^{-1}{\e^3}\Phi^2\al^{-4}\kappa^{-1}\,.
\ea
Inserting \eqref{wrefine14}, \eqref{wrefine15} and \eqref{wrefine16} into \eqref{wrefine13}, we conclude that, by choosing the small parameters appropriately and taking $\e$ sufficiently small
\bn
\|(u,v)\|^2_{L^2(\O)}\ls\al^{-6}\Phi^2\left(1+\kappa^{-2}\right)\e^3\,.
\en
This yields \eqref{wrefine0y} and completes the proof of Theorem \ref{thmexistence2x}.

\section*{Acknowledgments}

The authors thank Professor Jiaqi Yang and Dr. Jianfeng Zhao for helpful conversations. X. Pan is supported by the NSF of China under Grant No. 12471222.

\section*{Data availability statement}

Data sharing is not applicable to this article as no datasets were generated or analysed during the current study.

\section*{Conflict of interest statement}

The authors declare that they have no conflict of interest.

\end{document}